\documentclass[10pt,draft,twoside]{amsart}
\usepackage{amssymb}
\usepackage{latexsym}
\usepackage{amsfonts}
\usepackage{amsmath}
\usepackage{fancyheadings}
\usepackage{courier}

\newcommand{\K}{\mathcal K}

\newcommand{\E}{\mathcal E}

\newcommand{\cd}[2]{{\sf CD}_{#1}(#2)}

\newtheorem{theorem}{Theorem}[section]
\newtheorem{proposition}[theorem]{Proposition}
\newtheorem{lemma}[theorem]{Lemma}

\theoremstyle{definition}

\newtheorem{example}[theorem]{Example}

\newcommand{\alt}[1]{{\sf A}_{#1}}

\newcommand{\mat}[1]{{\sf M}_{#1}}
\newcommand{\sy}[1]{{\sf S}_{#1}}

\renewcommand{\sp}[2]{{\sf Sp}_{#1}(#2)}
\newcommand{\pomegap}[2]{{\sf P}\Omega^+_{#1}(#2)}

\newcommand{\sym}[1]{{\sf Sym}\,#1}
\newcommand{\supp}{\operatorname{Supp}}

\renewcommand{\wr}{\,{\sf wr}\,}

\newcommand{\aut}[1]{{\sf Aut}\,{#1}}

\newcommand{\inn}[1]{{\sf Inn}\,#1}

\newcommand{\cent}[2]{{\mathbb C}_{#1}(#2)}
\newcommand{\norm}[2]{{\mathbb N}_{#1}\left(#2\right)}

\newcommand{\psl}[2]{\mbox{\sf PSL}_{#1}(#2)}

\renewcommand{\sp}[2]{\mbox{\sf Sp}_{#1}(#2)}

\newcommand{\psu}[2]{\mbox{\sf PSU}_{#1}(#2)}

\newcommand{\Z}{\mathbb Z}

\renewcommand{\leq}{\leqslant}
\renewcommand{\geq}{\geqslant}

\begin{document}

\title[Intransitive Cartesian decompositions]{Intransitive Cartesian decompositions preserved by innately
transitive permutation groups}
\author{Robert W. Baddeley, Cheryl E. Praeger and Csaba Schneider}
\address[Baddeley]{Black Dog Software\\32 Arbury Road\\Cambridge CB4 2JE, UK}
\address[Praeger]{Department of Mathematics and Statistics\\
The University of Western Australia\\
35 Stirling Highway 6009 Crawley\\
Western Australia}
\address[Schneider]{Informatics Laboratory\\ 
Computer and Automation
Research Institute\\
The Hungarian Academy of Sciences\\
1111 Budapest L\'agym\'anyosi u.\ 11\\
Hungary}
\email{robert.baddeley@blackdogsoftware.co.uk, praeger@maths.uwa.edu.au,\hfill\break
csaba@maths.uwa.edu.au\protect{\newline} {\it WWW:}
www.maths.uwa.edu.au/$\sim$praeger, www.sztaki.hu/$\sim$schneider}

\begin{abstract}
A permutation group is innately transitive if it has a transitive
minimal normal subgroup, which is referred to as a plinth. 
We study the class of finite, innately
transitive permutation groups that can be embedded into wreath
products in product action. This investigation is carried out by
observing that such a wreath product preserves a natural Cartesian
decomposition 
of the underlying set. Previously we classified the possible
embeddings in the case where the innately transitive group projects onto a 
transitive subgroup of the top group. In this article we prove that 
innately transitive groups have at most three orbits
on an invariant  Cartesian decomposition. A consequence of this result is that if $G$ is an 
innately transitive subgroup of a wreath product in product action then
the natural projection of $G$ into the top group has at most two orbits.
\end{abstract}

\thanks{{\it Date:} 2 June 2004
The authors acknowledge the support of an Australian Research Council grant.}

\maketitle

\section{Introduction}\label{9.1}

The results presented in this 
paper play a key r\^ole in our research program to describe 
the Cartesian decompositions preserved by an innately
transitive permutation group. Recall that a permutation group $G$ is
said to be {\em innately transitive},
if $G$ has a transitive minimal normal subgroup $M$, which is
called a {\em plinth} of $G$. Innately transitive groups are
investigated in~\cite{bp}. 
The aim of our research is to describe
certain subgroups of wreath products in product action. We showed
in~\cite{recog} that these subgroups are best understood via studying
the natural Cartesian decomposition of the underlying set 
that is preserved by such a wreath product. In the same paper we
demonstrated the scope of 
this theory by describing transitive simple subgroups and their
normalisers in primitive wreath
products. Later in~\cite{transcs, 3types}
we described those innately transitive subgroups of wreath 
products in product action that project onto a transitive subgroup of the 
top group.

Here we consider the remaining case: 
we describe the innately transitive subgroups of wreath products in 
product action that project onto an intransitive subgroup of the top group. 
This amounts to saying that such a group acts intransitively on the 
corresponding Cartesian decomposition of the underlying set. 
The main result of this paper asserts that there are only two orbits in
this intransitive action.

\begin{theorem}\label{thA}
Suppose that $\Delta$ is a finite set, $|\Delta|\geq 2$, $\ell\geq 2$,
and $W$ is the wreath product $\sym\Delta\wr\sy\ell$ acting on
$\Delta^\ell$
in product action. Let $G$ be an innately transitive subgroup 
of $W$. Then the image of $G$ under the natural projection 
$W\rightarrow \sy\ell$ has at most two orbits on $\{1,\ldots,\ell\}$. 
\end{theorem}

The proof of Theorem~\ref{thA} is carried out by assuming that $G$
acts intransitively on the
underlying natural Cartesian decomposition of $\Delta^\ell$. Thus we
study Cartesian decompositions of sets that are acted upon
intransitively by an innately transitive permutation group. Though the
above-mentioned Cartesian decomposition of $\Delta^\ell$ is
homogeneous, that is, its elements have the same size, we do not
restrict our attention to this special case. Instead, we describe
innately transitive permutation groups acting intransitively on an
arbitrary Cartesian decomposition. The results of this study are
collected in Theorem~\ref{main}. Part~(iv) of Theorem~\ref{main}
implies Theorem~\ref{thA}, and also describes in more detail the
embedding in  Theorem~\ref{thA}.

The organisation of the paper is as follows. 
First in Section~\ref{sec2} we summarise those results of our 
previous work on Cartesian decompositions that will be used in this
paper. In the next section we
build the machinery that is necessary to investigate the scenario of
Theorem~\ref{thA}. Then we state
Theorem~\ref{main} which, as mentioned above, implies
Theorem~\ref{thA}. In order to prove our main theorem, we need results
about characteristically simple groups, and in Section~\ref{secfact}
we study the factorisations 
of such groups. In 
Section~\ref{secnorm} we prove several results about normalisers of
subgroups of characteristically simple groups.
Then in Sections~\ref{diagsec}, \ref{intsec}, and \ref{homsec} we treat Cartesian systems that are acted upon trivially
by a point stabiliser. Finally in Section~\ref{proof} we prove 
Theorem~\ref{main}.

Most of our results depend on the correctness of the finite simple
group classification. For instance, a lot of information on the
factorisations of simple and characteristically simple groups that depend
on this classification are used
throughout the paper.

The system of notation used in this paper is standard in permutation group
theory. Permutations act on the right: if $\pi$ is a permutation and
$\omega$ is a point then the image of $\omega$ under $\pi$ is denoted
$\omega\pi$.   If $G$ is a group acting on a
set $\Omega$ then $G^\Omega$ denotes the subgroup of $\sym\Omega$
induced by $G$. Further, if $\Gamma$ is a subset of $\Omega$ then $G_\Gamma$ 
and $G_{(\Gamma)}$ denote the setwise and the pointwise
stabilisers, respectively,  in $G$ of 
$\Gamma$. 

\section{Cartesian decompositions and Cartesian systems}\label{sec2}

A {\em Cartesian decomposition} of a set $\Omega$ is a set
$\{\Gamma_1,\ldots,\Gamma_\ell\}$ of proper partitions of $\Omega$ such that 
$$
|\gamma_1\cap\cdots\cap\gamma_\ell|=1\quad\mbox{for
all}\quad\gamma_1\in\Gamma_1,\ldots,\gamma_\ell\in\Gamma_\ell.
$$
This property implies that the following map is a well-defined bijection between $\Omega$ and
$\Gamma_1\times\cdots\times \Gamma_\ell$:
$$
\omega\mapsto(\gamma_1,\ldots,\gamma_\ell)\mbox{ where for }
i=1,\ldots,\ell,\ \gamma_i\in\Gamma_i\mbox{ is chosen so that }
\omega\in\gamma_i.
$$ 
Thus the set $\Omega$ can
naturally be
identified with the Cartesian product
$\Gamma_1\times\cdots\times\Gamma_\ell$. The number $\ell$ is called
the {\em index} of the Cartesian decomposition
$\{\Gamma_1,\ldots,\Gamma_\ell\}$. 


If $G$ is a permutation group acting on $\Omega$, then a Cartesian
decomposition $\E$ of $\Omega$ is said to be $G$-invariant, if the partitions in
$\E$ are permuted by $G$, and $\cd {}G$
denotes the set of $G$-invariant Cartesian decompositions of
$\Omega$. If $\E\in\cd {}G$ and $G$ acts on $\E$
transitively, then $\E$ is said to be a {\em transitive} $G$-invariant
Cartesian decomposition; otherwise it is called {\em intransitive}.  
The set of transitive $G$-invariant
Cartesian decompositions of $\Omega$ is denoted by $\cd{\rm
  tr}G$. 

The concept of a Cartesian decomposition was introduced by L.\ G.\ 
Kov\'acs in \cite{kov:decomp} where it was called a system of product imprimitivity.
Kov\'acs suggested that studying $\cd {\rm tr}G$ (using our terminology) was the appropriate way to identify 
wreath decompositions for finite primitive permutation groups $G$.
His papers
\cite{kov:decomp} and~\cite{kov:blowups} inspired our work.

Suppose that $G$ is an innately transitive
subgroup of $\sym\Omega$ with plinth $M$, and that $\E$ is a
$G$-invariant Cartesian decomposition of $\Omega$. 
In~\cite[Proposition~2.1]{recog} we proved that each of the $\Gamma_i$ is an
$M$-invariant partition of $\Omega$. 
Choose an element $\omega$ of $\Omega$ and let 
$\gamma_1\in\Gamma_1,\ldots,\gamma_\ell\in\Gamma_\ell$ be such that
$\{\omega\}=\gamma_1\cap\cdots\cap\gamma_\ell$; set $K_i=M_{\gamma_i}$. 
Then
\cite[Lemmas~2.2 and 2.3]{recog} imply that the set $\K_\omega(\E)=\{K_1,\ldots,K_
\ell\}$ is invariant
under conjugation by $G_\omega$, and, in addition,
\begin{eqnarray}\label{csdef1}
\bigcap_{i=1}^\ell K_i&=&M_\omega,\\
\label{csdef2}K_i\left(\bigcap_{j\neq
i}K_j\right)&=&M\quad \mbox{for all}\quad i\in\{1,\ldots,\ell\}.
\end{eqnarray}

For an arbitrary transitive permutation group $M$ on $\Omega$ and a point $\omega\in\Omega$, a set $\K=\{K_1,\ldots,K_\ell\}$ of proper subgroups
of $M$ is said to be 
a {\em Cartesian system of subgroups with respect to $\omega$} 
for $M$, if~\eqref{csdef1} and \eqref{csdef2}~hold.

\begin{theorem}[Theorem~1.4 and Lemma~2.3~\cite{recog}]\label{bij}
Let $G\leq \sym\Omega$ be an innately transitive permutation group with
plinth $M$. For a fixed $\omega\in\Omega$ 
the correspondence $\E\mapsto\K_\omega(\E)$ is a bijection between the
set of $G$--invariant Cartesian decompositions of $\Omega$
and the set of $G_\omega$--invariant
Cartesian systems of subgroups for $M$ with respect to
$\omega$. Moreover the $G_\omega$--actions on $\E$ and on $\K_\omega(\E)$
are equivalent.
\end{theorem}

With $G\leq\sym\Omega$, $M$, and $\omega\in\Omega$ as above, let $\K$ be a
$G_\omega$-invariant  Cartesian system of subgroups for $M$ with
respect to $\omega$. Then
Theorem~\ref{bij} implies that $\K=\K_\omega(\E)$ for some 
$G$-invariant Cartesian decomposition $\E$ of $\Omega$. In fact, $\E$ 
consists of the $M$-invariant partitions 
$\{(\omega^K)^m\ |\ m\in M\}$
where $K$ runs through the elements of $\K$. 
This Cartesian decomposition
is usually denoted $\E(\K)$.

Using this
theory we were able to describe in~\cite{recog} 
those innately transitive subgroups 
of wreath products that have a simple plinth. This led to
a classification of transitive simple subgroups of
wreath products in product action
(see~\cite[Theorem~1.1]{recog}). Then in~\cite{transcs,3types} we
extended this classification and described innately transitive
subgroups of such wreath products that project onto a transitive
subgroup of the top group.

Suppose now that $M=T_1\times\cdots\times T_k$ where the $T_i$ are
groups, and $k\geq1$. For
$I\subseteq\{T_1,\ldots,T_k\}$, $\sigma_I:M\rightarrow\prod_{T_i\in
I}T_i$ denotes the natural projection map. We also write
$\sigma_i$ for $\sigma_{\{T_i\}}$. 
A subgroup $X$ of $M$ is said to be a {\em strip} if, for each
$i=1,\ldots,k$, either $\sigma_i(X)=1$ or $\sigma_i(X)\cong X$. The set
of all $T_i$ such that $\sigma_i(X)\neq 1$ is called the {\em support} of
$X$ and is denoted $\supp X$, and $|\supp X|$ is called the {\em length} of
$X$. If $T_m\in\supp X$ then we also say that
$X$ {\em covers} $T_m$. Two strips $X_1$ and $X_2$ are said to be 
{\em disjoint} if $\supp X_1\cap\supp X_2=\emptyset$. 
A strip $X$ is said to be {\em full} if
$\sigma_i(X)=T_i$ for all $T_i\in\supp X$,
and it is called {\em non-trivial} if
$|\supp X|\geq 2$. A subgroup $K$ of $M$ is said to be {\em subdirect} with 
respect to the direct decomposition $T_1\times\cdots\times T_k$ if
$\sigma_i(K)=T_i$ for all $i$. If $M$ is a finite, non-abelian,
characteristically simple group,
then a subgroup $K$ is said to be {\em subdirect} if it is subdirect with
respect to the finest direct decomposition of $M$ (that is, the
product decomposition with simple groups as factors).

The importance of strips is highlighted by the following result, which
is usually referred to as Scott's Lemma (see the appendix of~\cite{scott}).

\begin{lemma}\label{scott}
Let $M$ be a direct product of finitely many non-abelian, finite
simple groups and $H$ a
subdirect subgroup of $M$. Then $H$ is a direct product of pairwise
disjoint full strips of $M$.
\end{lemma}

Let $M=T_1\times\cdots\times T_k$ be a finite, non-abelian, 
characteristically
simple group, where $T_1,\ldots,T_k$ are the simple normal subgroups of
$M$, each isomorphic to the same simple group $T$. If $K$ is a subgroup of 
$M$ and $X$ is a strip in $M$ such that $K=X\times\sigma_{\{T_1,\ldots,T_k\}
\setminus\supp X}(K)$ then we say that $X$ is {\em involved} in $K$. 
A strip $X$ is said to be involved in a Cartesian system $\K$ for $M$
if $X$ is involved in some element of $\K$. 
Note that in this case~\eqref{csdef2} 
implies that $X$ is involved in a unique element of $\K$.

A non-abelian plinth of an innately transitive group 
$G$ has the form $M=T_1\times\cdots\times T_k$ where
the $T_i$ are  finite, non-abelian, simple groups. Let $\E\in\cd{}G$ and let
$\K_\omega(\E)$ be a corresponding Cartesian system
$\{K_1,\ldots,K_\ell\}$ for $M$ with respect to $\omega$. Then equation~\eqref{csdef2} implies that, 
for all $i\leq k$ and $j\leq\ell$,
\begin{equation}\label{simpfact}
\sigma_i(K_j)\left(\bigcap_{j'\neq j}\sigma_i(K_{j'})\right)=T_i.
\end{equation}
In particular this means that if $\sigma_i(K_j)$ is a proper subgroup
of $T_i$ then $\sigma_i(K_{j'})\neq\sigma_i(K_j)$ for all
$j'\in\{1,\ldots,\ell\}\setminus\{j\}$. 
It is thus important to understand the following sets of subgroups:
\begin{equation}\label{f}
\mathcal F_i(\E,M,\omega)=\{\sigma_i(K_j)\ |\ j=1,\ldots,\ell,\ \sigma_i(K_j)
\neq T_i\}.
\end{equation}
From our remarks above, $|\mathcal F_i(\E,M,\omega)|$ is the number of 
indices $j$ such that
$\sigma_i(K_j)\neq T_i$. 
The set $\mathcal F_i(\E,M,\omega)$ is independent of $i$ up to isomorphism, 
in the sense that if $i_1,\ i_2\in\{1,\ldots,k\}$ and $g\in G_\omega$ 
are such that
$T_{i_1}^g=T_{i_2}$ then $\mathcal
F_{i_1}(\E,M,\omega)^g=\{L^g\ |\ L\in\mathcal F_{i_1}(\E,M,\omega)\}=
\mathcal F_{i_2}(\E,M,\omega)$. This argument
also shows that the subgroups in $\mathcal F_{i_1}(\E,M,\omega)$ are actually 
$G_\omega$-conjugate to the subgroups in $\mathcal F_{i_2}(\E,M,\omega)$.

The set $\cd{\rm tr}G$ is further subdivided according to
the structure of the subgroups in the corresponding Cartesian systems
as follows. The sets $\mathcal F_i=\mathcal F_i(\E,M,\omega)$ are defined
in~\eqref{f}.
\begin{eqnarray*}
\cd{\rm S}G&=&\{\E\in\cd{\rm tr}G\ |\ \mbox{the elements of $\K_\omega(\E)$
are subdirect subgroups in $M$}\};\\
\cd 1G&=&\{\E\in\cd{\rm tr}G\ |\ |\mathcal F_i|=1\mbox{ and $\K_\omega(\E)$
involves no non-trivial, full strip}\};\\
\cd {\rm 1S}G&=&\{\E\in\cd{\rm tr}G\ |\ |\mathcal F_i|=1\mbox{ and $\K_\omega(\E)$
involves non-trivial, full strips}\};\\
\cd {2\sim}G&=&\{\E\in\cd{\rm tr}G\ |\ |\mathcal F_i|=2\mbox{ and
the $\mathcal F_i$ contain two $G_\omega$-conjugate subgroups}\};\\
\cd{2\not\sim}G&=&\{\E\in\cd{\rm tr}G\ |\ |\mathcal F_i|=2\mbox{ and
the subgroups in $\mathcal F_i$ are not $G_\omega$-conjugate}\};\\
\cd 3G&=&\{\E\in\cd{\rm tr}G\ |\ |\mathcal F_i|=3\}.
\end{eqnarray*}

At first glance, it seems that the definitions of the classes $\cd {\rm S}G$, 
$\cd 1G$, $\cd{\rm 1S}G$, $\cd {2\sim}G$, $\cd {2\not\sim}G$, and $\cd 3G$ 
may depend on the choice of the Cartesian system, and hence on the choice 
of the point $\omega$. However, the
following result, proved in~\cite[Theorems~6.2 and~6.3]{transcs}, 
shows that this is not the case, and also implies that these classes
form a partition of $\cd{\rm tr}G$. A finite permutation group is said to be
{\em quasiprimitive} if all its non-trivial normal subgroups are
transitive. We also say that a quasiprimitive group has compound
diagonal type, if it has a unique minimal normal subgroup, which is
non-abelian, and in which a point stabiliser is a non-simple 
subdirect subgroup.

\begin{theorem}[6-class Theorem]\label{5class}
If $G$ is a finite, innately transitive permutation group with a non-abelian
plinth $M$, then the classes $\cd 1G$, $\cd{\rm S}G$, $\cd {\rm 1S}G$,
$\cd{2\sim}G$, $\cd {2\not\sim}G$, and $\cd 3G$
are independent of the choice of the point $\omega$ used in their definition.
They form a partition of $\cd{\rm tr}G$, and moreover, if $M$ is simple,
then $\cd{\rm tr}G=\cd {2\sim}G$. Suppose, in addition, that $T$ is
the common isomorphism type of the simple direct factors of $M$. Then 
the following all hold.
\begin{enumerate}
\item[(a)] 
The group $G$ is quasiprimitive of compound diagonal type if and only if $\cd{\rm S}G\neq\emptyset$.
\item[(b)] If $\cd {\rm 1S}G\cup\cd {2\sim}G\neq\emptyset$,
then $T$ 
has a factorisation with two isomorphic, proper subgroups and 
is isomorphic to one of the groups $\alt 6$, $\mat{12}$, 
$\pomegap 8q$, or $\sp 4{2^a}$ with $a\geq2$. If $\E\in \cd {\rm 1S}G\cup\cd {2\sim}G$ then the subgroups in
$\mathcal F_i(\E,M,\omega)$ are isomorphic to the groups $A$ and $B$ in the corresponding 
line of Table~$\ref{isomfact}$. 
\item[(c)] If $\cd {2\not\sim}G\neq\emptyset$, then $T$ admits a
factorisation $T=AB$ with $A, B$ proper subgroups.
\item[(d)] If $\cd 3G\neq\emptyset$, then $T$ is isomorphic to one of the 
groups $\sp{4a}2$ with $a\geq 2$, $\pomegap 83$, or $\sp 62$. If $\E\in\cd 3G$ then the 
subgroups in $\mathcal F_i(\E,M,\omega)$ are isomorphic to the groups $A$, $B$, and $C$ in the
corresponding line of  Table~$\ref{tablemult}$.
\end{enumerate}
\end{theorem}

\section{Intransitive Cartesian decompositions}\label{mainthsec}

In this section we state Theorem~\ref{main}, which can be
viewed as a qualitative characterisation of innately transitive
groups acting intransitively on a Cartesian decomposition. In
particular Theorem~\ref{thA} follows from the first assertion and
part~(iv) of this result. Before we can state this theorem, we
introduce some notation which will also be used in later parts of this paper.

Suppose that $G$ is an innately transitive permutation group acting on 
$\Omega$ with plinth $M$ and let 
$\E=\{\Gamma_1,\ldots,\Gamma_\ell\}$ be a $G$-invariant Cartesian 
decomposition of $\Omega$ on which $G$ acts intransitively. 
It follows from~\cite[Proposition~5.1]{prae:inc} that $M$ is non-abelian. 
Suppose that 
$\Xi_1,\ldots,\Xi_s$ are the $G$-orbits on $\E$, and that
$\K_\omega(\E)=\{L_1,\ldots,L_\ell\}$ 
is the Cartesian system of subgroups for $M$ with respect to
some fixed $\omega\in\Omega$. 

For 
$i=1,\ldots,s$ set
$$
K_i=\bigcap_{\Gamma_j\in \Xi_i}L_j.
$$
For partitions $A_1,\ldots,A_d$ of a set $\Omega$, the 
infimum $\inf\{A_1,\ldots,A_d\}$ of these partitions is defined as the 
partition
$$
\inf\{A_1,\ldots,A_d\}=\left\{\gamma_1\cap\cdots\cap\gamma_d\ 
|\ \gamma_1\in A_1,\ldots,\gamma_d\in A_d\right\}.
$$
The proof that $\inf \{A_1,\ldots,A_d\}$ is a partition of
$\Omega$ is easy and is left to the reader. 
We note that $\inf\{A_1,\ldots,A_d\}$ is
the coarsest partition that refines each of the $A_i$, and hence is
the infimum of $A_1,\ldots,A_d$ with respect to the natural partial
order on the set of partitions of $\Omega$.  
Let $\Omega_i$ denote $\inf\Xi_i$ for $i=1,\ldots,s$ and let
$\bar\E=\{\Omega_1,\ldots,\Omega_s\}$. If 
$\gamma\in\Gamma_j$ for some $\Gamma_j\in\Xi_i$
then $\gamma$ is a union of blocks from $\Omega_i$ and we set 
$$
\bar\gamma=\{\delta\ |\ \delta\in\Omega_i,\  
\delta\subseteq\gamma\},\quad \bar\Gamma_j=\{\bar\gamma\
|\ \gamma\in\Gamma_j\}\quad\mbox{and}\quad
\bar\Xi_i=\{\bar\Gamma_j\ |\ \Gamma_j\in\Xi_i\}.
$$ 

It turns out, as shown in the following theorem, that $\bar\E$ is a
Cartesian decomposition of $\Omega$ acted upon trivially by
$G$. Further, each $G$-invariant partition $\Omega_i$ in $\bar\E$ admits a
$G$-invariant, transitive Cartesian decomposition, namely
$\bar\Xi_i$. Thus the study of the original intransitive decomposition
$\E$ can be carried out via the study of a $G$-trivial decomposition, and the study of several transitive Cartesian
decompositions. 
This idea is made more explicit in Theorem~\ref{main}.
The concepts full factorisation, full strip factorisation, and
strong multiple factorisation occurring in the statement 
are defined in Section~\ref{secfact}. 

\begin{theorem}\label{main}
The number $s$ of $G$-orbits on $\E$ is at most~$3$. The partitions
$\Omega_i$ are $G$-invariant and the set $\bar\E=\{\Omega_1,\ldots,\Omega_s\}$
is a Cartesian decomposition of $\Omega$ on which $G$ acts trivially.
For $i=1,\ldots,s$, the subgroup $K_i$ is the stabiliser in $M$ of the 
block in $\Omega_i$ containing $\omega$. Moreover, for $i=1,\ldots,s$,
the $M$-action on
$\Omega_i$ is faithful and $\bar\Xi_i\in\cd{\rm tr}{G^{\Omega_i}}$. 
\begin{enumerate}
\item[(i)] If $\bar\Xi_i\in\cd{\rm S}{G^{\Omega_i}}$, for
some $i\in\{1,\ldots,s\}$, then $s=2$. 
Further, if, say, 
$\bar\Xi_{1}\in \cd{\rm S}{G^{\Omega_1}}$,
then
$(M,K_1,K_2)$ is a full strip factorisation, and $\bar\Xi_{2}\in\cd 1{G^{\Omega_{2}}}$. 
\item[(ii)] If $\bar\Xi_i\in\cd{2\not\sim}{G^{\Omega_i}}$ for some
$i\in\{1,\ldots,s\}$, then $s=2$, and, for 
all $j\in\{1,\ldots,k\}$, the group $T_j$ 
and the subgroups of $\mathcal F_j(\E,M,\omega)$ 
are as in Table~$\ref{tablemult}$. If  $\bar\Xi_1\in\cd{2\not\sim}{G^{\Omega_1}}$ then
$\bar\Xi_2\in\cd{1}{G^{\Omega_1}}$.
\item[(iii)] We have, for all $i\in\{1,\ldots,s\}$, that
$\bar\Xi_i\not\in\cd {\rm 1S}{G^{\Omega_i}}\cup \cd{2\sim}{G^{\Omega_i}}\cup
\cd {3}{G^{\Omega_i}}$.
\item[(iv)] If $\E$ is homogeneous then $s=2$, $\bar\Xi_i\in\cd 1{G^{\Omega_i}}$ 
  for $i=1,\ 2$, and $(M,\{K_1,K_2\})$ is a full factorisation.
\item[(v)] If $s=3$ then $\bar\Xi_i\in\cd 1{G^{\Omega_i}}$ 
for $i=1,\ 2,\ 3$, and $(M,\{K_1,K_2,K_3\})$ is a strong multiple factorisation.
\end{enumerate}
\end{theorem}

The general part of Theorem~\ref{main} follows from the following
result. The rest of the assertions made in Theorem~\ref{main} will be
verified in  Section~\ref{proof}.

\begin{proposition}\label{1ststep}
Let $G$, $M$, $\Omega$, $\omega$, 
$\E$, $\bar\Xi_1,\ldots,\bar\Xi_s$, $\Omega_1,\ldots,\Omega_s$, 
and $K_1,\ldots,K_s$ 
be as above. 
Then $\bar\E=\{\Omega_1,\ldots,\Omega_s\}$ 
is a $G$-invariant Cartesian decomposition of $\Omega$
such that the Cartesian system $\K_\omega(\bar\E)$ coincides with
$\{K_1,\ldots,K_s\}$. Moreover, for each $i$, 
$\Omega_i$ is a $G$-invariant partition of $\Omega$, $K_i$ 
is normalised by $G_\omega$, $\bar\Xi_i\in\cd{\rm
  tr}{G^{\Omega_i}}$,  and $M$ acts faithfully on $\Omega_i$.  
\end{proposition}
\begin{proof}
In the first two paragraphs we prove  that 
$\Omega_1$ is a $G$-invariant partition of
$\Omega$, that $K_1$ is the stabiliser in $M$ of the part of
$\Omega_1$ containing $\omega$ (and hence $G_\omega$ normalises $K_1$),
and that $M$ is faithful on $\Omega_1$. 
The proofs for the other $\Omega_i$ are identical. 
Suppose that
$\Xi_1=\{\Gamma_1,\ldots,\Gamma_m\}$ and let
$\gamma_1\in\Gamma_1,\ldots,\gamma_m\in\Gamma_m$ be the blocks
containing $\omega$. Then, by the definition of $\Omega_1$,
$\gamma_1\cap\cdots\cap\gamma_m\in\Omega_1$. 
It follows from the definition of the infimum that $\Omega_1$ is a
partition of $\Omega$. Let $g\in G$. 
Then
$g$ permutes $\Gamma_1,\ldots,\Gamma_m$ among themselves, and so
$\{\gamma_1^g,\ldots,\gamma_m^g\}\cap\Gamma_i$ is a singleton for each
$i\in\{1,\ldots,m\}$. Thus 
$(\gamma_1\cap\cdots\cap\gamma_m)^g=\gamma_1^g\cap\cdots\cap\gamma_m^g\in\Omega_1$,
and so $\Omega_1$ is $G$-invariant.

Next we prove that $K_1$ is
the stabiliser in $M$ of the element $\gamma_1\cap\cdots\cap\gamma_m$
in $\Omega_1$. Now
$K_1=L_1\cap\cdots\cap L_m$ where 
$L_j$ is the stabiliser in $M$ of $\gamma_j$.
Hence $K_1$ stabilises
$\gamma_1\cap\cdots\cap\gamma_m$, the block in $\Omega_1$ that
contains $\omega$. Now suppose that some element $g\in M$ stabilises
$\gamma_1\cap\cdots\cap\gamma_m$. The definition of a Cartesian system
implies that $\gamma_1\cap\cdots\cap\gamma_m$ is non-empty. As
$\gamma_1,\ldots,\gamma_m$ are blocks of imprimitivity 
for the $M$-actions on $\Gamma_1,\ldots,\Gamma_m$, respectively, it
follows that $g$ fixes each of $\gamma_1,\ldots,\gamma_m$ setwise. 
Thus $g\in
L_1\cap\cdots\cap L_m$. Therefore $K_1$ is the
stabiliser in $M$ of $\gamma_1\cap\cdots\cap\gamma_m$. As $\Omega_1$ is a $G$-invariant
partition of $\Omega$,  
$K_1$ is normalised
by $G_\omega$. 
Moreover, since
$K_1\neq M$ it follows that $\bigcap_{g\in G}K_1^g$ is a normal
subgroup of $G$ properly contained in $M$. As $M$ is a minimal normal
subgroup of $G$, this implies that $\bigcap_{g\in G}K_1^g=1$, and so
$M$ acts faithfully on $\Omega_1$.

We now claim that $\bar\E=\{\Omega_1,\ldots,\Omega_s\}$ is a Cartesian
decomposition of $\Omega$.
Let $\delta_1\in\Omega_1,\ldots,\delta_s\in\Omega_s$. Because
of the definition of the $\Omega_i$, there are
$\gamma_1\in\Gamma_1,\ldots,\gamma_\ell\in\Gamma_\ell$ such that
$\delta_1\cap\cdots\cap\delta_s=\gamma_1\cap\cdots\cap\gamma_\ell$.
As the $\Gamma_i$ form a Cartesian decomposition of $\Omega$, we
obtain that $|\gamma_1\cap\cdots\cap\gamma_\ell|=1$, and so $
|\delta_1\cap\cdots\cap\delta_s|=1$. Thus $\bar\E$ is a Cartesian
decomposition of $\Omega$. Since each of the $\Omega_i$ is a
$G$-invariant partition of $\Omega$, the group $G$
acts trivially on $\bar\E$. Since $K_i$ is the stabiliser in $M$ of
the part of $\Omega_i$ containing $\omega$, it follows that
$\{K_1,\ldots,K_m\}$ is the Cartesian system $\K_\omega(\bar\E)$.

Finally we prove that
$\bar\Xi_i\in\cd{\rm tr}{G^{\Omega_i}}$ for each $i=1,\ldots,s$, and,
as usual, we show this for $i=1$.  Recall that
$\Xi_1=\{\Gamma_1,\ldots,\Gamma_m\}$. 
Suppose that $\bar\Gamma_1,\ldots,\bar\Gamma_m$ are as above, 
and let $\bar\gamma_1\in\bar\Gamma_1,\ldots,
\bar\gamma_m\in\bar\Gamma_m$ corresponding to the elements
$\gamma_1\in\Gamma_1,\ldots,\gamma_m\in\Gamma_m$, respectively. 
Since $\gamma_1\cap\cdots\cap\gamma_m$ is a block in $\Omega_1$, we have
$|\bar\gamma_1\cap\cdots\cap\bar\gamma_m|=1$. This shows that $\bar\Xi_1$ 
is a Cartesian decomposition of $\Omega_1$. The 
$G$-actions on $\Xi_1$ and $\bar \Xi_1$ are naturally equivalent, and, as $G$ is 
transitive on $\Xi_1$, we obtain that $\bar\Xi_1\in\cd{\rm tr}{G^{\Omega_1}}$. 
\end{proof}

\section{Factorisations of simple and characteristically simple
groups}\label{secfact}

The factorisations of simple and characteristically simple groups play an
important r\^ole in this paper, especially in the proof of 
Theorem~\ref{main}. Such factorisations were studied 
earlier in~\cite{bad:fact,charfact}. In this section we summarise and extend 
the results proved in these papers.

A {\em group factorisation} is a pair $(G,\{A,B\})$ where $G$ is a group and
$A,\ B$ are subgroups of $G$ such that $AB=G$. 
In this situation we also say that $\{A,\ B\}$ is a
factorisation of $G$, and we often write that $G=AB$ is a
factorisation. A factorisation 
is called
{\em non-trivial} if both $A$ and $B$ are proper subgroups.
In this paper we only consider non-trivial
factorisations.

Let $M=T_1\times\cdots\times T_k$ be a finite, non-abelian,
characteristically simple group where $T_1,\ldots,T_k$ are pairwise
isomorphic, 
simple normal subgroups. Then a  factorisation $M=K_1K_2$ is said to
be a {\em full factorisation} if, for each $i\in\{1,\ldots,k\}$,
\begin{enumerate}
\item[(a)] the subgroups $\sigma_i(K_1),\ \sigma_i(K_2)$ are proper
subgroups of $T_i$;
\item[(b)] the orders $|\sigma_i(K_1)|$, $|\sigma_i(K_2)|$,
  and $|T_i|$ are divisible by the same set of primes.
\end{enumerate}

Full factorisations of simple and characteristically simple
groups were classified in~\cite{bad:fact} and~\cite{charfact},
respectively. The following result is a short summary of what we
need to know about such factorisations to prove the results  in this paper.
 
\begin{theorem}\label{fullth}
Suppose that $k\geq 1$ and 
$T_1,\ldots,T_k$ are pairwise isomorphic, finite, non-abelian
simple groups, and set $M=T_1\times\cdots\times T_k$. If $M=K_1K_2$ is
a full factorisation then
$$
\sigma_1(K_j)'\times\cdots\times\sigma_k(K_j)'\leq K_j\quad\mbox{for}\quad j\in\{1,\ 2\}.
$$
Further, for each $i\in\{1,\ldots,k\}$, the pair
$(T_i,\{\sigma_i(K_1), \sigma_i(K_2)\})$ occurs as $(T,\{A,B\})$
in one of the lines of 
Table~$\ref{full_simple}$. 
\end{theorem}

\begin{center}
\begin{table}[ht]
$$
\begin{array}{|l|c|c|c|}
\hline
 & T & A & B\\
\hline
1 & \alt 6&\alt 5&\alt 5 \\
\hline
2 & \mat{12}&\mat{11}&\mat{11},\ \psl 2{11} 
 \\
\hline
3 & \pomegap 8q,\ q\geq 3&\Omega_7(q)&\Omega_7(q)\\
\hline
4 & \pomegap 82&\sp 62&\alt 7,\ \alt 8,\ \sy 7,\ \sy 8,\ \sp
 62,\ \Z_2^6\rtimes\alt 7,\ \Z_2^6\rtimes\alt 8\\
\cline{3-4}
&&\alt 9 & \alt 8,\ \sy 8,\ \sp 62,\  \Z_2^6\rtimes\alt 7,\ \alt 8,\ \Z_2^6\rtimes\alt 8\\
\hline
5 & \sp 4q,\ q\geq 4\mbox{ even}
& \sp 2{q^2}.2 & \sp2{q^2}.2,\ \sp 2{q^2} \\
\hline
\end{array}
$$
\caption{Full factorisations $\{A,B\}$ of finite simple groups $T$}
\label{full_simple}
\end{table}
\end{center}

An important subfamily of full factorisations consists of the
factorisations of non-abelian, finite simple groups with two
isomorphic subgroups. We will use the extra details about these
factorisations given below.

\begin{lemma}
\label{factlem}
Let $T$ be a finite simple group and $A,\ B$ proper subgroups of $T$ such
that $|A|=|B|$ and $T=AB$. Then the following hold.
\begin{enumerate}
\item[(i)] The isomorphism types of 
$T$, $A$, and $B$ are as in
Table~$\ref{isomfact}$, and $A$, $B$ are isomorphic, 
maximal subgroups of $T$.
\item[(ii)] There is an automorphism $\vartheta\in\aut T$ that
interchanges $A$ and $B$.
\item[(iii)] We have $\norm T{A'\cap B'}=\norm T{A\cap B}=A\cap B$ and
$\cent T{A'\cap B'}=\cent T{A\cap B}=1$.
\end{enumerate}
\end{lemma}
\begin{center}
\begin{table}[ht]
$$
\begin{array}{|l|c|c|}
\hline
 & T & A,\ B\ \\
\hline
1 & \alt 6&\alt 5\\
\hline
2 & \mat{12}&\mat{11}
 \\
\hline
3 & \pomegap 8q &\Omega_7(q)\\
\hline
4 & \sp 4q,\ q\geq 4\mbox{ even}
& \sp 2{q^2}.2  \\
\hline
\end{array}
$$
\caption{Factorisations of finite simple groups in Lemma~\ref{factlem}}
\label{isomfact}
\end{table}
\end{center}
\begin{proof}
Parts~(i) and~(ii) were proved in~\cite[Lemma~5.2]{recog}, 
and the same result implies that $A\cap B$ is self-normalising in 
$T$. It is not hard to see that the proof of~\cite[Lemma~5.2]{recog} can 
be used, after minor alteration,  to verify that $\norm T{A'\cap B'}=A\cap B$.
In particular $\cent T{A\cap B}={\sf Z}(A\cap
B)$ and $\cent T{A\cap B}=\cent{A\cap B}{A'\cap B'}$. 
We obtain from the~\cite{atlas} in rows~1--2, from~\cite[3.1.1(vi)]{kleidman}
in row~3, and from~\cite[3.2.1(d)]{lps:max} in row~4 of
the table that $A\cap B$ is a centerless
group and, in row~4, that $\cent{A\cap B}{A'\cap B'}=1$.
\end{proof}

Let $M=T_1\times\cdots\times
T_k$ be a finite, non-abelian, 
characteristically simple group as above. 
For subgroups  $K_1,\ldots,K_\ell$ of $M$,
the pair $(M,\{K_1,\ldots,K_\ell\})$ is said to be a 
{\sl strong multiple factorisation} if, 
for all $i\in\{1,\ldots,k\}$ and all pairwise distinct 
$j_1,\ j_2,\ j_3\in\{1,\ldots,\ell\}$,
\begin{enumerate}
\item[(a)] $\sigma_i(K_1),\ldots,\sigma_i(K_\ell)$ are proper
  subgroups of $T_i$; and
\item[(b)] $K_{j_1}(K_{j_2}\cap K_{j_3})=K_{j_2}(K_{j_1}\cap
  K_{j_3})=K_{j_3}(K_{j_1}\cap K_{j_2})=M$.
\end{enumerate}

The following theorem, combining~\cite[Table~V]{bad:fact}
and~\cite[Theorem~1.7, Corollary~1.8]{charfact}, gives a characterisation of strong multiple
factorisations of characteristically simple groups.

\begin{theorem}\label{th2}
A strong multiple factorisation of a finite characteristically simple
group contains exactly three subgroups.
If $M$ is a non-abelian, characteristically simple group with simple
normal subgroups $T_1,\ldots,T_k$, and  
$(M,\{K_1,K_2,K_3\})$ is a strong multiple factorisation, then
$\sigma_1(K_i)'\times\cdots\times\sigma_k(K_i)'\leq
K_i$ for $i=1,\ 2,\ 3$,
and, for $i=1,\ldots,k$, the pair
$(T_i, \{\sigma_i(K_1), \sigma_i(K_2), \sigma_i(K_3)\})$ occurs
as $(T, \{A, B, C\})$ in one of the lines of 
Table~$\ref{tablemult}$. Further, if one of the lines $1$--$2$ of
Table~$\ref{tablemult}$ is valid then 
$\sigma_1(K_i)\times\cdots\times\sigma_k(K_i)=
K_i$ for $i=1,\ 2,\ 3$.
\end{theorem}

\begin{center}
 \begin{table} 
$$
 \begin{array}{|c|c|c|c|c|}
 \hline
 & T & A & B & C \\
 \hline
 1& \sp {4a}2,\ a\geq 2 &\sp {2a}4.2 & {\sf O}^-_{4a}(2) & {\sf
 O}^+_{4a}(2)\\
 \hline
 2 & \pomegap 83 & \Omega_7(3) & \Z_3^6\rtimes \psl 43 & \pomegap 82\\
 \hline
 3 & \sp 62 & {\sf G}_2(2) & {\sf O}^-_{6}(2) & {\sf O}^+_{6}(2)\\
 \cline{3-5}
  && {\sf G}_2(2)' & {\sf O}^-_{6}(2) & {\sf O}^+_{6}(2)\\
 \cline{3-5}
  && {\sf G}_2(2) & {\sf O}^-_{6}(2)' & {\sf O}^+_{6}(2)\\
 \cline{3-5}
  && {\sf G}_2(2) & {\sf O}^-_{6}(2) & {\sf O}^+_{6}(2)'\\
 \hline
 \end{array}
 $$
 \caption{Strong multiple factorisations $\{A,B,C\}$ of finite simple
 groups $T$}\label{tablemult}
 \end{table}
\end{center}

The concept of a full strip factorisation is defined for the purposes of this
paper.
For the characteristically simple group
$M=T_1\times\cdots\times T_k$ and proper subgroups $D$, $K$, the
triple $(M,D,K)$ is said to
be a {\em full strip factorisation} if 
\begin{enumerate}
\item[(i)] $M=DK$;
\item[(ii)] $D$ is a direct product of pairwise disjoint, non-trivial full
strips;
\item[(iii)]  for all $i,\ j\in\{1,\ldots,k\}$, $\sigma_i(K)$ is a proper subgroup of $T_i$ and
$\sigma_i(K)\cong\sigma_j(K)$.
\end{enumerate}

The following lemma is shows that            in a full strip 
factorisation each full strip has length~2.

\begin{lemma}\label{fsf}
If $(M,D,K)$ is a full strip factorisation of a finite,
characteristically simple group $M$, then each non-trivial, full
strip involved in $D$ has length $2$.
\end{lemma}
\begin{proof}
Suppose without loss of generality that $X$ is a non-trivial full
strip involved in $D$, covering $T_1,\ldots,T_s$ for some $s\geq
2$. We let $I=\{T_1,\ldots,T_s\}$. Then
$\sigma_{I}(D)=X$ and the
factorisation  $X\sigma_{I}(K)=T_1\times\cdots\times T_s$
holds. Then \cite[Lemma~4.3]{transcs} implies that $s\leq 3$, and if
$s=3$ then the simple direct factor $T$ of $M$ 
admits a strong multiple factorisation involving three subgroups 
isomorphic to the subgroups $\sigma_i(K)$, for $i=1,\ 2,\ 3$. 
On the other hand,
Table~\ref{tablemult} shows that finite simple groups do not admit
strong multiple factorisations with isomorphic subgroups. This is a
contradiction, and hence $s=2$.
\end{proof}

The next result provides the link between the concepts of a full
factorisation and a full strip 
factorisation.

\begin{theorem}[Theorem~1.5~\cite{charfact}]\label{diag}
 Let $M=T_1\times\cdots\times T_{2k}$ be a characteristically
 simple group, where the $T_i$ are non-abelian, simple groups, 
$\varphi_i:T_i\rightarrow T_{i+k}$ an isomorphism for
 $i=1,\ldots,k$, and set
 $$
 D=\{(t_1,\ldots,t_{k},\varphi_1(t_{1}),\ldots,\varphi_k(t_{k}))\ |\
 t_1\in T_1,\ldots,t_{k}\in
 T_k\}.
 $$
 If $(M,D,K)$ is a full strip
 factorisation, then
 $(T_i,\{\sigma_i(K),\varphi_i^{-1}(\sigma_{i+k}(K))\})$ is a factorisation
 of $T_i$ with isomorphic subgroups for all $i\in\{1,\ldots,k\}$, and $\prod_{i=1}^{2k}\sigma_i(K)'\leq K$.
 \end{theorem}

The following 
useful result from~\cite{bad:quasi} shows that in a non-trivial
factorisation of a non-abelian characteristically simple group, it is
not possible for both factors to be direct products of pairwise
disjoint strips.

\begin{lemma}[{\cite[Lemma~2.2]{bad:quasi}}]\label{2.1}
Suppose that $M=T_1\times \cdots\times T_k$ is a direct product of
isomorphic non-abelian, simple groups $T_1,\ldots,T_k$. Suppose that
$A_1,\ldots,A_{m}$ are non-trivial pairwise disjoint strips in $M$,
and so are $B_1,\ldots,B_n$. Then $M\neq (A_1\times\cdots\times
A_m)(B_1\times\cdots\times B_n)$. 
\end{lemma}

\section{Normalisers in direct products}\label{secnorm}

In this section we collect together some facts about normalisers 
of subgroups in direct products that will be used in our analysis. 
It is easy to see that the
normaliser in a direct product $G_1\times\cdots\times G_k$ of a
subgroup $H$ is contained in
$\norm{G_1}{\sigma_1(H)}\times\cdots\times\norm{G_k}{\sigma_k(H)}$. 
Moreover, if $H$ is the direct product
$\sigma_1(H)\times\cdots\times\sigma_k(H)$ then 
$\norm{G_1\times\cdots\times
G_k}H=\norm{G_1}{\sigma_1(H)}\times\cdots\times\norm{G_k}{\sigma_k(H)}$. The following simple lemma extends this observation to a more
general situation.

\begin{lemma}\label{samenorm}
Let $G_1,\ldots,G_k$ be groups, set $G=G_1\times\cdots\times G_k$
and for $i=1,\ldots,k$ let $H_i$ be a subgroup of $G_i$. Let $H$ be a
subgroup of $G$ such that $H_1\times\cdots\times H_k\lhd H$,
the factor $\norm G{H_1\times\cdots\times H_k}/(H_1\times\cdots\times H_k)$ is abelian, and $\norm
{G_i}{\sigma_i(H)}=\norm{G_i}{H_i}$. Then 
$$
\norm GH=\norm
G{H_1\times\cdots\times H_k}=\norm{G_1}{H_1}\times\cdots\times\norm{G_k}{H_k}.
$$ 
\end{lemma}
\begin{proof}
As $\norm
{G_i}{\sigma_i(H)}=\norm{G_i}{H_i}$, it follows from the remarks above
that
$$
\norm GH\leq \prod_{i=1}^k\norm{G_i}{H_i}=\norm
G{H_1\times\cdots\times H_k}.
$$
On the other hand, as
$H_1\times\cdots\times H_k\lhd H$ and 
$\norm G{H_1\times\cdots\times H_k}/(H_1\times\cdots\times
H_k)$ is abelian, $\norm
G{H_1\times\cdots\times H_k}\leq\norm GH$. Therefore equality holds.
\end{proof}

The following lemma, proved in~\cite[Lemma~3.5]{3types}, determines 
the normaliser of a strip.

\begin{lemma}\label{normstriplemma}
Let $G_1,\ldots, G_k$ be isomorphic groups, $\varphi_i:G_1\rightarrow G_i$
an isomorphism for $i=2,\ldots,k$, $H_1$ a subgroup of $G_1$, and  $H=\{(h,\varphi_2(h),\ldots,\varphi_k(h))\ |\ t\in H_1\}$
a non-trivial
strip in $G_1\times\cdots\times G_k$. Then 
$$
\norm {G_1\times\cdots\times G_k}H=\left\{(t,c_2\varphi_2(t),\ldots,c_k\varphi_k(t))\ |\ t\in\norm {G_1}{H_1},\ c_i\in\cent
{G_i}{\varphi_i(H_1)}\right\}.
$$
\end{lemma}

We use the results above to derive some facts concerning normalisers of
the subgroups that occur in Table~\ref{groupsinttable}.

\begin{proposition}\label{normcor1}
Suppose that  $M=T_1\times\cdots\times T_k\cong T^k$ 
is a characteristically simple group and
$(M,\{K_1,K_2\})$ is a
full factorisation such that, for all $i$, the pair
$(T_i,\{\sigma_i(K_1),\sigma_i(K_2)\})$ is as $(T,\{A,B\})$ in one of
the rows of Table~$\ref{groupsinttable}$.
\begin{itemize}
\item[(a)] If $T$ is as in one of rows $1$--$3$ of Table~$\ref{groupsinttable}$
then $K_1$, $K_2$, and $K_1\cap K_2$ 
are self-normalising in $M$.
\item[(b)] If row~$4$ of Table~$\ref{groupsinttable}$ is valid
then, for $j=1,\ 2$, 
 we have $\norm M{K_j}=\prod_i\sigma_i(K_j)$ and
$\norm M{K_1\cap K_2}=\norm M{K_1}\cap\norm M{K_2}$.
\end{itemize}
\end{proposition}
\begin{proof}
(a) In this case the $\sigma_i(K_j)$ are perfect and, 
by Theorem~\ref{fullth},
$K_j=\prod_i\sigma_i(K_j)$ for $j=1,\ 2$, and Table~\ref{groupsinttable}
shows that
$\sigma_i(K_j)$ is self normalising in $T_i$ for all
$i\in\{1,\ldots,k\}$ and $j\in\{1,\ 2\}$.
Therefore $K_1$ and $K_2$ are self-normalising.
Further,
$$
K_1\cap
K_2=\prod_{i=1}^k\sigma_i(K_1\cap
K_2)=\prod_{i=1}^k\sigma_i(K_1)\cap\sigma_i(K_2).
$$ 
Using Lemma~\ref{factlem} and the Atlas~\cite{atlas}, we obtain that $\norm
{T_i}{\sigma_i(K_1)\cap\sigma_i(K_2)}=\sigma_i(K_1)\cap\sigma_i(K_2)$ for
all $i$. Thus
\begin{multline*}
\norm M{K_1\cap K_2}=\norm
M{\prod_{i=1}^k\sigma_i(K_1)\cap\sigma_i(K_2)}=\prod_{i=1}^k\norm
{T_i}{\sigma_i(K_1)\cap\sigma_i(K_2)}\\=\prod_{i=1}^k\sigma_i(K_1)\cap\sigma_i(K_2)=\prod_{i=1}^k
\sigma_i(K_1\cap K_2)=K_1\cap K_2.
\end{multline*}

(b) Now assume that $T\cong\sp 4q$ for some $q\geq 4$ even and let
$j\in\{1,2\}$.  By Theorem~\ref{fullth}, $K_j'=\prod_i\sigma_i(K_j)'$, and 
we can read off
from Table~\ref{groupsinttable} that
$\norm{T_i}{\sigma_i(K_j)'}=\norm{T_i}{\sigma_i(K_j)}=\sigma_i(K_j)$
for all $i\in\{1,\ldots,k\}$. 
As $\norm M{K_j'}/K_j'$ is elementary abelian and $\norm
M{K_j'}\geq K_j\geq K_j'$, Lemma~\ref{samenorm} gives
$\norm M{K_j'}= \norm M{K_j}$.
On the other hand,
$$
\norm
M{K_j'}=\prod_{i=1}^k\norm {T_i}{\sigma_i(K_j)'}=\prod_{i=1}^k\norm {T_i}{\sigma_i(K_j)}=\prod_{i=1}^k\sigma_i(K_j).
$$

Now Theorem~\ref{fullth} shows that $K_1'\cap
K_2'=\prod_i\sigma_i(K_1'\cap K_2')$. We also obtain from
Theorem~\ref{fullth} and 
Lemma~\ref{factlem} that 
$$
\norm{T_i}{\sigma_i(K_1\cap
K_2)}=\norm{T_i}{\sigma_i(K_1'\cap
K_2')}=\norm{T_i}{\sigma_i(K_1)}\cap\norm{T_i}{\sigma_i(K_2)}. 
$$
Thus
\begin{multline*}
\norm M{K_1'\cap
K_2'}=\prod_{i=1}^k\norm{T_i}{\sigma_i(K_1'\cap
K_2')}
=\prod_{i=1}^k\left(\norm
{T_i}{\sigma_i(K_1)}\cap\norm {T_i}{\sigma_i(K_2)}\right)\\
=\prod_{i=1}^k\norm
{T_i}{\sigma_i(K_1)}\cap\prod_{i=1}^k\norm {T_i}{\sigma_i(K_2)}=\norm M{K_1}\cap\norm M{K_2}.
\end{multline*}
As $\norm M{K_1'\cap K_2'}/(K_1'\cap K_2')$ is abelian,
and $K_1'\cap K_2'\leq K_1\cap K_2\leq \norm M{K_1'\cap K_2'}$,
Lemma~\ref{samenorm} implies that $\norm M{K_1\cap K_2}=\norm M{K_1'\cap
K_2'}=\norm M{K_1}\cap \norm M{K_2}$. 
\end{proof}

\begin{proposition}\label{normcor2}
Let $M=T_1\times\cdots\times T_{2k}=T^{2k}$, $D$ and $K$ be as in Theorem~$\ref{diag}$, and suppose that
$DK=M$. 
\begin{itemize}
\item[(a)] If $T$ is as in one of the rows~$1$--$3$ of
Table~$\ref{isomfact}$
then 
$K$ and $K\cap D$ are self-normalising in $M$. 
\item[(b)] If $T$ is as in row~$4$ of Table~$\ref{isomfact}$ then
$\norm MK=\prod_i\sigma_i(K)$ and 
$\norm M{K\cap D}=D\cap\norm MK$.
\end{itemize}
\end{proposition}
\begin{proof}
(a) Theorem~\ref{diag} implies that $K$ is the direct product of its
projections onto the $T_i$, and by Table~\ref{isomfact}  these projections are
self-normalising in $T$. Hence $K$
is self-normalising. Suppose that
$X=\{(t,\varphi_1(t))\ |\ t\in T_1\}$ is a full strip involved in $D$
where $\varphi_1:T_1\rightarrow T_{k+1}$ is an isomorphism. 
Then 
$\sigma_{\{T_1,T_{k+1}\}}(K\cap D)=\{(t,\varphi_1(t))\ |\
t\in\sigma_1(K)\cap\varphi_1^{-1}(\sigma_{k+1}(K))\}$. By Theorem~\ref{diag},
$\sigma_1(K)\varphi_1^{-1}(\sigma_{k+1}(K))=T_1$ is a factorisation with
isomorphic subgroups.
Hence Lemma~\ref{factlem} implies that
$\sigma_1(K)\cap\varphi_1^{-1}(\sigma_{k+1}(K))$ is self-normalising in
$T_1$ 
and that $\cent{T_1}{\sigma_1(K)\cap\varphi_1^{-1}(\sigma_{k+1}(K))}=1$.
Thus Lemma~\ref{normstriplemma} yields that $\sigma_{\{T_1,T_{k+1}\}}(K\cap D)$ is
self-normalising in $T_1\times T_{k+1}$. Similar argument shows that $\sigma_{\{T_i,T_{i+k}\}}(K\cap D)$ is
self-normalising in $T_i\times T_{i+k}$ for all $i\in\{1,\ldots,k\}$. 
As $K\cap D=\sigma_{\{T_1,T_{k+1}\}}(K\cap D)\times\cdots\times\sigma_{\{T_k,T_{2k}\}}(K\cap
D)$ 
we obtain that $K\cap D$ is also
self-normalising in $M$.

(b) The argument which was used in part~(b) of
Proposition~\ref{normcor1} to compute the normalisers of 
$K_1$ and $K_2$ can be used to verify the claim about $\norm MK$. Using the argument in part~(a), 
it is easy to check that $\norm M{D\cap K'}=D\cap \norm MK$. 
Since 
$\norm M{D\cap K'}/(D\cap K')$ is
abelian, $\norm M{D\cap K'}\geq D\cap K\geq D\cap K'$, and $\norm{T_i\times
T_{i+k}}{\sigma_{\{T_i,T_{i+k}\}}(D\cap K)}=\norm{T_i\times
T_{i+k}}{\sigma_{\{T_i,T_{i+k}\}}(D\cap K')}$, we obtain from Lemma~\ref{samenorm}
that $\norm M{D\cap K'}=\norm M{D\cap K}$. 
\end{proof}

\section{Cartesian systems involving non-trivial
strips}\label{diagsec}

We use the notation introduced in Section~\ref{sec2}. 
Let us start with a motivating example.

\begin{example}\label{stex}
Let $T$ be a finite simple group with two proper, isomorphic subgroups $A$
and $B$, such that $T=AB$. The possibilities for $T$, $A$, and $B$ are
in Table~\ref{isomfact}. Suppose that $k$ is even and set $M=T^k$. 
Let $\tau$ be an element of $\aut T$ 
interchanging $A$ and $B$; such a $\tau$ exists by
Lemma~\ref{factlem}.
Consider the following two subgroups
of $M$:
$$
\bar K_1=\{(t_1,t_1,\ldots,t_{k/2},t_{k/2})\ |\
t_1,\ldots,t_{k/2}\in T\}\quad\mbox{and}\quad \bar K_2=(A\times B)^{k/2}.
$$
We obtain from Theorem~\ref{diag} that $M=K_1K_2$. 
Identify $M$ with $\inn M$ in $\aut M\cong\aut T\wr\sy k$, and set 
$$
\bar G=M\left(\norm{\aut M}{\bar K_1}\cap\norm{\aut M}{\bar
K_2}\right). 
$$
We
claim that $M$ is a minimal normal subgroup of $\bar G$, or
equivalently, $\bar G$ is transitive by conjugation on the simple direct factors
of $M$. Note that $\sigma_1(\bar K_2)=A$ and $\sigma_{2}(\bar K_2)=B$, 
and the
automorphism 
$(\tau,\tau,1,\ldots,1)(1,2)$
of $M$ interchanges the first two simple direct factors of $M$, while normalising $\bar K_1$
and $\bar K_2$. 
Also the automorphism $(1,3,\ldots,k-1)(2,4,\ldots,k)$ of $M$ cyclically
permutes the blocks determined by the strips in $\bar K_1$, and
normalises both $\bar K_1$ and $\bar K_2$. Therefore the subgroup of $\aut M$
generated by these two outer automorphisms is transitive
on the simple direct factors of $M$, and, in addition, normalises
$\bar K_1$ and $\bar K_2$.  Hence $M$ is a
minimal normal subgroup of $\bar G$ and, since $\cent{\aut M}M=1$, it
is the unique such minimal normal subgroup.

Set $G_0=\norm{\aut M}{\bar K_1}\cap\norm{\aut M}{\bar K_2}$ so that
$MG_0=\bar G$. As $\bar K_1$
and $\bar K_2$ are self-normalising in $M$, we obtain  $M\cap G_0=\bar
K_1\cap \bar K_2$. 
Therefore, by~\cite[Lemma~4.1]{3types}, the $M$-action on the coset space
$[M:\bar K_1\cap \bar K_2]$ can be extended to $\bar G$ with point stabiliser
$G_0$.
Moreover, $\bar K_1$ and $\bar K_2$ form a Cartesian system
for $M$ acted upon trivially by $G_0$. Consequently this action of $\bar
G$ preserves an intransitive $\bar G$-invariant Cartesian decomposition,
such that one of the subgroups, namely $\bar K_1$, in the corresponding Cartesian system
$\{\bar K_1,\bar K_2\}$ is the direct product of disjoint strips.
\end{example}

Our aim in this section is to describe the intransitive, pointwise
$\bar G$-invariant Cartesian
decompositions whose Cartesian systems involve a non-trivial full
strip. If $K_i$ involves such a strip for some $i$ 
then $\sigma_j(K_i)=T_j$ for some $j\in\{1,\ldots,k\}$. Without loss
of generality we may suppose that $1=i=j$. In this case we obtain the
following theorem.

\begin{theorem}\label{intrstripd2}
Let $G$, $M$, and $\K$ be as in Section~$\ref{mainthsec}$, and assume that
$\sigma_1(K_1)=T_1$.
Then $s=2$ and $(M,K_1,K_2)$ is a full strip
factorisation. 
In particular, the isomorphism types of
$T$ and $\sigma_i(K_2)$ are as in Table~$\ref{isomfact}$. 
Further,
$K_2'=\sigma_1(K_2)'\times\cdots\times\sigma_k(K_2)'$, and if $T$ is
not as in row~$4$ of Table~$\ref{isomfact}$ then
$K_2=\sigma_1(K_2)\times\cdots\times\sigma_k(K_2)$. 
\end{theorem}
\begin{proof}
Since $G_\omega$ normalises $K_1$ and acts
transitively on the $T_i$, 
we have that $\sigma_i(K_1)=T_i$ for all $i$. If $T_i\leq
K_1$ for some $i$ then, for the same reason, 
$T_i\leq K_1$ for all $i$ and so $K_1=M$,
which is impossible. Hence, by Scott's Lemma~\ref{scott}, 
$K_1$ is a direct product
of non-trivial full strips. If $\sigma_i(K_j)=T_i$ for some $i$ and
some $j\neq
1$ then the same argument shows that $K_j$ is also a direct product of
non-trivial, full strips. However, Lemma~\ref{2.1} implies that $K_1K_j\neq
M$, which violates the defining properties of Cartesian systems. Thus
$\sigma_i(K_j)$ is a proper subgroup of $T_i$ for all $i$ and all $j\geq
2$. 

Since $G_\omega$ normalises $K_2$,
$\sigma_i(K_2)\cong\sigma_j(K_2)$ for all $i$ and $j$. Thus
$(M,K_1,K_2)$ is a full strip factorisation.
Lemma~\ref{fsf} implies that all strips involved in $K_1$ have length 2. 

We now show that $s=2$. Suppose on the contrary that $s\geq
3$. Let $X$ be a strip in $K_1$ whose support is, without loss of generality,
$\{T_1,T_2\}$. Then 
$X=\{(t,\alpha(t))\ |\ t\in T_1\}$ 
for some isomorphism $\alpha:T_1\rightarrow T_2$.
For $i\geq 2$, $\sigma_1(K_i)$ and
$\sigma_2(K_i)$ are proper subgroups of $T_1$ and
$T_2$, respectively, and it
follows from Theorem~\ref{diag}
that $(T_1,\{\sigma_1(K_i),\alpha^{-1}(\sigma_2(K_i))\})$ 
is a factorisation
with isomorphic subgroups. 
As $K_2$ and $K_3$ are normalised by $G_\omega$, so is
their intersection $K_2\cap K_3$. 
Hence $\sigma_1(K_2\cap
K_3)\cong\sigma_2(K_2\cap K_3)$. Since
$K_1(K_2\cap K_3)=M$ and $\sigma_{\{1,2\}}(K_1)$ is the full strip $X$,
we obtain from~\cite[Lemma~2.1]{charfact} that 
$(T_1,\{\sigma_1(K_2\cap K_3),\alpha^{-1}(\sigma_2(K_2\cap K_3))\})$ is
also a full factorisation with isomorphic subgroups. In such
factorisations the subgroups involved are maximal subgroups of $T_1$ (see
Table~\ref{isomfact}), and so $\sigma_1(K_2\cap K_3)$ and
$\alpha^{-1}(\sigma_2(K_2\cap K_3))$ are maximal subgroups of $T_1$. However,
$\sigma_1(K_2\cap K_3)\leq \sigma_1(K_2)\cap \sigma_1(K_3)$, which, as
$\sigma_1(K_2)$ and $\sigma_1(K_3)$ are proper subgroups of $T_1$,
implies that $\sigma_1(K_2\cap
K_3)$, $\sigma_1(K_2)$, and $\sigma_1(K_3)$ coincide. Hence
$\sigma_1(K_2K_3)=\sigma_1(K_2)\sigma_1(K_3)<T_1$ which is a
contradiction, as $K_2K_3=M$. Thus $s=2$. The rest of the theorem
follows from Theorem~\ref{diag} and from the fact that the subgroups
$A$ and $B$ in rows 1--3 of Table~\ref{isomfact} are perfect.
\end{proof}

If $K_1$ is a subdirect subgroup of $M$, then 
we prove that $\cent GM$ is small, in fact, in most cases $\cent
GM=1$ and $G$ is quasiprimitive. If $G$ is a permutation group with a
unique minimal normal subgroup $M$, such that $M$ is transitive, then
$G$ is quasiprimitive. Moreover if $M$ is not simple, a
point stabiliser in $M$ is non-trivial and is not a subdirect subgroup
of $M$, then $G$ is said to have quasiprimitive type {\sc Pa};
see~\cite{bad:quasi}.

\begin{proposition}\label{soccentdiag}
Let $G$, $M$, and $\K$ be as in Section~$\ref{mainthsec}$. 
Assume that $\sigma_1(K_1)=T_1$. 
If the group $T$ is as in rows~$1$--$3$ of Table~{\rm \ref{isomfact}}, then
$\cent{\sym\Omega}M=1$, and in
particular $G$ is quasiprimitive of type {\sc Pa}. If $T$ is as in row~$4$ then $\norm
M{M_\omega}=K_1\cap\norm M{K_2}$, and 
$$
\cent{\sym\Omega}M\cong (K_1\cap\norm M{K_2})/(K_1\cap K_2)\cong \norm
M{K_2}/K_2.
$$
In particular $\cent{\sym\Omega}M$ is an
elementary abelian $2$-group of rank at most $k/2$, and all minimal normal subgroups of $G$
different from $M$ are elementary abelian $2$-groups.
\end{proposition}
\begin{proof}
By Theorem~\ref{intrstripd2}, $s=2$, and so, $M_\omega=K_1\cap K_2$. 
Note that by~\cite[Theorem~4.2A]{dm} 
$$
\cent{\sym\Omega}{M}\cong\norm M{M_\omega}/M_{\omega}=\norm M{K_1\cap
K_2}/(K_1\cap K_2).
$$
If $T$ is as in one of the rows~1--3 of Table~\ref{isomfact}, then
Proposition~\ref{normcor2} implies that $K_1\cap K_2$ is
self-normalising in $M$, and hence $\cent{\sym\Omega}M=1$. This
implies that $M$ is the unique minimal normal subgroup in $G$, and so
$G$ is quasiprimitive. As $K_1$ involves a non-trivial full strip,
$k\geq 2$. Moreover, it follows from Table~\ref{isomfact} that
$M_\omega\neq 1$ and $M_\omega$ is not a subdirect subgroup of
$M$. Thus $G$ has quasiprimitive type {\sc Pa}. If $T$ is as in
row~5, then, again by Proposition~\ref{normcor2}, we only have to prove that
$(K_1\cap\norm M{K_2})/(K_1\cap K_2)$ and $\norm M{K_2}/K_2$ are
isomorphic.

Recall that
$K_1K_2=M$,
and so $\norm M{K_2} =(K_1\cap \norm M{K_2})K_2$. By the second
isomorphism theorem, $\norm M{K_2}/K_2\cong (K_1\cap \norm
M{K_2})/(K_1\cap K_2)$ under the isomorphism $\psi(xK_2)=x'(K_1\cap
K_2)$ where $x=x'x''\in\norm M{K_2}$ with $x'\in K_1\cap \norm
M{K_2}$, $x''\in K_2$. A proof that $\psi$ is well-defined and
is an isomorphism can be found in most group theory textbooks; see for
instance~\cite[3.12 Satz]{huppert}. 
\end{proof}

\section{Bounding the number of orbits in an intransitive Cartesian system}\label{intsec}

We apply the results of the last section to prove that $s\leq
3$.

\begin{theorem}\label{3th}
The index $s$ of the Cartesian system $\K$ in Theorem~$\ref{main}$ is at most $3$. Further, if $s=3$ then $(M,\{K_1,K_2,K_3\})$ is a strong
multiple factorisation. Hence, in this case, for all $i$, the
factorisation $(T_i,\{\sigma_i(K_1),\sigma_i(K_2),\sigma_i(K_3)\})$ is 
in Table~$\ref{tablemult}$. Moreover if $T$ is as in the first two rows then 
$$
K_i=\sigma_1(K_i)\times\cdots\times\sigma_k(K_i)\quad\mbox{for}\quad
i=1,\ 2,\ 3,
$$
while if $T$ is as in the third row then
$$
\sigma_1(K_i)'\times\cdots\times\sigma_k(K_i)'\leq K_i\leq \sigma_1(K_i)\times\cdots\times\sigma_k(K_i)\quad\mbox{for}\quad
i=1,\ 2,\ 3.
$$
\end{theorem}
\begin{proof}
If $\sigma_i(K_j)=T_i$ for some $i$ and $j$
then, by Theorem~\ref{intrstripd2}, we  
have $s=2$.  Therefore we may assume without loss of generality that all projections $\sigma_i(K_j)$ are proper in $T_i$. 
Then $K_1,\ldots,K_s$ form a strong multiple factorisation of $M$.
Thus, by Theorem~\ref{th2}, $s\leq 3$. 

If $s=3$, then, by Theorem~\ref{intrstripd2}, $\sigma_i(K_j)<T_i$ for
all $i$ and $j$. Thus, by~\eqref{csdef2}, 
$\{K_1,K_2,K_3\}$ is a strong multiple
factorisation of $M$, and $T_i$,
$\sigma_i(K_1)$, $\sigma_i(K_2)$, $\sigma_i(K_3)$ are as in
Table~\ref{tablemult}. The assertions about the $K_i$ follow
from Theorem~\ref{th2}. 
\end{proof}

A generic example with $s=3$ can easily be constructed as follows.

\begin{example}\label{exsm}
Let ${A,B,C}$ be maximal subgroups of a finite simple group $T$
forming a strong multiple factorisation of $T$, and let $\bar K_1=A^k$,
$\bar K_2=B^k$, $\bar K_3=C^k$ be the corresponding subgroups of $M=T^k$.
Then $(T,\{A,B,C\})$ and $(M,\{\bar K_1,\bar K_2,\bar K_3\})$ are strong multiple
factorisations. 
Identify $M$ with $\inn M$ in
$\aut M$, and let 
$$
\bar G=M\left(\norm {\aut M}{\bar K_1}\cap\norm {\aut M}{\bar
K_2}\cap\norm {\aut M}{\bar K_3}\right).
$$
Since the cyclic subgroup of $\aut M$ generated by the 
automorphism $$
\tau:(x_1,\ldots,x_k)\mapsto(x_k,x_1,\ldots,x_{k-1})
$$ 
is transitive on the simple direct factors of $M$ and
normalises $\bar K_1$, $\bar K_2$, and $\bar K_3$, we have that $M$ is a minimal normal subgroup of $\bar G$. Moreover,
since $\cent{\aut M}M=1$, $M$ is the unique minimal normal
subgroup of $\bar G$. 

If $G_0=\norm{\aut M}{\bar K_1}\cap\norm{\aut M}{\bar
K_2}\cap\norm{\aut M}{\bar K_3}$ then
$MG_0=\bar G$ and, since $\bar K_1$,
$\bar K_2$, and $\bar K_3$ are self-normalising in $M$, $M\cap G_0=\bar
K_1\cap \bar K_2\cap\bar K_3$. 
Therefore, by~\cite[Lemma~4.1]{3types}, the $M$-action on the coset space
$[M:\bar K_1\cap \bar K_2\cap\bar K_3]$ can be extended to $\bar G$ with point stabiliser
$G_0$.
Moreover, $\{\bar K_1,\bar K_2,\bar K_3\}$ is a Cartesian system
for $M$ acted upon trivially by $G_0$. Consequently this action of $\bar
G$ preserves an intransitive $\bar G$-invariant Cartesian decomposition
given by the Cartesian system
$\{\bar K_1,\bar K_2,\bar K_3\}$.
\end{example}

The defining properties of $\K$ give us some useful constraints on
$T$. For instance if the $K_i$ involve no non-trivial, full strips,
then $T_i=\sigma_i(K_j)\sigma_i(K_m)$ for all $i,\ j,\ m$
such that $j\neq m$. In particular $T$ has a proper factorisation, and
so, for example, $T\not\cong\psu{2d+1}q$ unless $(d,q)\in\{(1,3),\
(1,5),\ (4,2)\}$. Many sporadic simple groups can also be
excluded. See the tables in~\cite{lps:max}.

In general it is difficult to give a complete description of Cartesian
decompositions that involve no strips. However we can give such a
description when the initial intransitive Cartesian decomposition $\E$
is homogeneous, This 
leads to the proof of Theorem~\ref{thA}.
Describing the remaining case would be more difficult than
finding all factorisations of finite simple groups, as demonstrated by
the following generic example.

\begin{example}
Let $T$ be a finite simple group, $k\geq 1$, and set $M=T^k$.
Let $\{A,B\}$ be a non-trivial factorisation of the group $T$ and set $K_1=A^k$,
$K_2=B^k$. Then clearly $K_1K_2=T^k$, and the base group $T^k$ is the unique minimal normal
subgroup of $G=T\wr \sy k$.  Consider the 
coset action of $G$ on $\Omega=[G:G_0]$
where $G_0=(A\cap B)\wr \sy k$. Then $K_1\cap K_2=(A\cap B)^k=M\cap
G_0$, and $K_1,\ K_2$ are normalised
by $G_0$, so they  give rise to a $G$-invariant 
intransitive Cartesian decomposition of $\Omega$ with index 2.
\end{example}

The example above shows that a detailed description of all Cartesian
decompositions preserved by an innately transitive group would first
require determining all factorisations of finite simple
groups. But even assuming that such a classification is available,
determining the relevant factorisations of characteristically simple
groups is still  a difficult task. In the cases that we investigate in
the remainder of this paper the required factorisations of the $T_i$ were readily
available. The subgroups of these factorisations were almost simple or
perfect which made possible an explicit description of the occurring
factorisations of $M$.

\section{Intransitive homogeneous Cartesian decompositions}\label{homsec}

The aim of this section is to describe homogeneous, intransitive
Cartesian decompositions preserved by an innately transitive group.
Such Cartesian
decompositions need to be studied if we want to investigate embeddings
of innately transitive groups in wreath products in product action. 
First we note, using the notation introduced for Theorem~\ref{main}, 
that $|\Gamma_i|=|\Gamma_j|$ for all $i,\ j$. 
Then, for each
$i\in\{1,\ldots,\ell\}$, $m=|\Gamma_i|$ (independent of $i$), and there
is an integer $\ell_i$ such that 
$|\Omega_i|=|M:K_i|=|\Gamma_1|^{\ell_i}=m^{\ell_i}$ for all $i\in\{1,\ldots,s\}$. 

\begin{theorem}\label{homth}
Let $G$, $M$, $\E$, and $\K$ be as in Section~$\ref{mainthsec}$.
Assuming that $\E$ is homogeneous, we have
$\sigma_i(K_j)<T_i$ for all $i$ and $j$. Further, in this case, $s=2$ and
$(M,\{K_1,K_2\})$ is a full factorisation.
\end{theorem}
\begin{proof}
Let us first prove that $\sigma_i(K_j)<T_i$ for all $i$ and $j$. Suppose
without loss of generality that $\sigma_1(K_1)=T_1$. Then Theorem~\ref{intrstripd2}
implies that $s=2$, $K_1$ is the direct product of strips of length 2,
and $K_2'=\sigma_1(K_2)'\times\cdots\times\sigma_k(K_2)'\leq
K_2$. Recall that there exist non-negative integers $m,\ \ell_1,\ \ell_2$ such
that $[M:K_1]=m^{\ell_1}$ and $[M:K_2]=m^{\ell_2}$. Since
$|K_1|\cong|T|^{k/2}$ we have $[M:K_1]=|T|^{k/2}$, and so all primes that
divide $|T|$ will also divide $m$. Since
$K_2'\leq K_2$ and  $K_2'$ is the direct product of its projections $\sigma_i(K_2)'$,
it follows that $|M:K_2'|=|T_1:\sigma_1(K_2)'|^k$, and so all prime
divisors $p$ of $|T|$ divide $[T_1:\sigma_1(K_2)']$. This is not the case if
$T\cong \alt 6$ or $T\cong\mat {12}$ (take $p=5$ in both cases). If
$T\cong \pomegap 8q$ and $\sigma_1(K_2)\cong \Omega_7(q)$ then 
$$
|T|=\frac{1}{d^2}q^{12}(q^6-1)(q^4-1)^2(q^2-1)\quad\mbox{and}\quad
|T_1:\sigma_1(K_2)'|=\frac{1}{d}q^3(q^4-1)
$$
where $d=(2,q-1)$. By Zsigmondy's Theorem (see~\cite[\S 2.4]{lps:max}), there exists a prime $r$
dividing $q^6-1$ and not dividing $q^4-1$, whence $r$
divides $|T|$ but not $|T_1:\sigma_1(K_2)'|$. When $T\cong\sp 4q$ with
$q$ even, $q\geq 4$, then $\sigma_1(K_2)'\cong\sp 2{q^2}$, so
$$
|T|=q^4(q^4-1)(q^2-1)\quad\mbox{and}\quad
 |T_1:\sigma_1(K_2)'|=q^2(q^2-1).
$$
Using Zsigmondy's theorem we find that $q^4-1$ has a prime divisor $r$
that does not divide $q^2-1$. Thus $r$ divides $|T|$ but not
$|T_1:\sigma_1(K_2)'|$.  Therefore $\sigma_i(K_j)<T_i$ for all $i$ and
$j$.

For all distinct $i,\ j\in\{1,\ldots,s\}$ we have $M=K_iK_j$, and
hence $m^{\ell_i}=|M:K_i|=|K_j:K_i\cap K_j|$ divides $K_j$. It follows
that every prime divisor of $m$ divides $|K_j|$. Let $p$ be a prime
divisor or $|T|$. Then $p$ divides $|M|$. Since $|M:K_j|=m^{\ell_j}$,
either $p$ divides $|K_j|$ or $p$ divides $m$, and in the latter case we
also obtain that $p$ divides $|K_j|$. By Proposition~\ref{1ststep}, $G_\omega$
normalises $K_j$, and since $G=MG_\omega$, $G_\omega$ acts
transitively by conjugation on $\{T_1,\ldots,T_k\}$. It follows that,
for $1\leq i\leq k$, 
the projections $\sigma_i(K_j)$ are pairwise isomorphic, proper
subgroups of $T_i$. Thus, since $p$ divides $|K_j|$, we deduce that $p$ divides
$\sigma_i(K_j)$, for each $i$. Hence each prime divisor of $|T|$
divides $|\sigma_i(K_j)|$ for all $i$ and $j$. Set $Q_j=\sigma_1(K_j)$
for $j=1,\ldots,s$. 

If $s\geq 3$ then, since $\K$ is a Cartesian system, 
$(T_1,\{Q_1,\ldots,Q_s\})$ is a strong multiple
factorisation (see the paragraph before Theorem~\ref{th2}). 
Moreover, since $|T|$, $|Q_i|$, $|Q_j|$ are divisible by the same
primes, $(T_1,\{Q_i,Q_j\})$ is a full factorisation for
all $i\neq j$. Comparing Tables~\ref{full_simple} and~\ref{tablemult}, 
we find that no strong
multiple factorisation of a finite simple group exists in which any two of the subgroups form a full
factorisation. 
Hence we obtain that $s=2$ and $(M,\{K_1,K_2\})$ is a full
factorisation.
\end{proof}

\begin{center}
\begin{table}[ht]
\[
\begin{array}{|l|c|c|c|}
\hline
 & T & A & B \\
\hline
1 & \alt 6&\alt 5&\alt 5\\
\hline
2 & \mat{12}&\mat{11}& \mat{11},\ \psl 2{11}\\
\hline
3 & \pomegap 8 q&\Omega_7(q)&\Omega_7(q)\\
\hline
4& \sp 4q,\ q\geq 4\ \hbox{and $q$ even}&\sp
4{q^2}.2&\sp 2{q^2}.2\\
\hline
\end{array}
\]
\caption{The table for Theorem~\ref{groups-int}}
\label{groupsinttable}
\end{table}
\end{center}

\begin{theorem}\label{groups-int}
Let $G$, $M$, $T_1,\ldots,T_k$, $\E$, 
and $\K$ be as in Section~$\ref{mainthsec}$. If $\E$ is homogeneous,
then, for all $i\in\{1,\ldots,k\}$, 
$(T_i,\{\sigma_i(K_1),\sigma_i(K_2)\})$ is a factorisation 
$(T,\{A,B\})$ as in one of the lines of 
Table~{\rm \ref{groupsinttable}}. If $T$ is as in rows~$1$--$3$,
then $K_i=\sigma_1(K_i)\times \cdots\times\sigma_k(K_i)$ for $i=1,\ 2$.
Moreover in this case $\cent{\sym\Omega}M=1$, and in
particular $G$ is quasiprimitive. If $T$ is as in row~$4$ then 
\begin{equation}\label{commnorm}
\sigma_1(K_i)'\times\cdots\times\sigma_k(K_i)'\leq K_i\leq
\sigma_1(K_i)\times\cdots\times\sigma_k(K_i)=\norm
M{K_i}\quad\mbox{for}\quad i=1,\ 2,
\end{equation}
and
$$
\cent{\sym\Omega}M\cong (\norm M{K_1}\cap\norm M{K_2})/(K_1\cap K_2).
$$
In particular $\cent{\sym\Omega}M$ is an
elementary abelian $2$-group of rank at most $k$, and all minimal normal subgroups of $G$
different from $M$ are elementary abelian $2$-groups.
\end{theorem}
\begin{proof}
Set $A=\sigma_1(K_1)$ and $B=\sigma_1(K_2)$, so that, by Theorem~\ref{homth}, $T_1=AB$ is a full factorisation.
We have to eliminate all full factorisations of $T$ which are not
contained in Table~\ref{groupsinttable}. These involve the group
$T=\sp 4q$ or $\pomegap 82$, 
and we consider these families separately.

Suppose first that $T\cong\sp 4q$ with $q$ even, $q\geq 4$. If $A$ and
$B$ are isomorphic to $\sp 2{q^2}.2$ then line~4 of
Table~\ref{groupsinttable} is valid. Suppose that $A$,
say,  is
isomorphic to $\sp 2{q^2}$.  
Then $K_1$ is isomorphic to $A^k\cong (\sp 2{q^2})^k$. 
As the factorisation $K_1K_2=M$ holds, we must have, for all $i$,  that
$\sigma_i(K_2)\cong \sp 2{q^2}\cdot 2$, and hence $|K_1|<|K_2|$. 
For a positive integer $n$ and a prime $p$ let $n_p$ denote the
exponent of the largest $p$-power dividing $n$.  Recall that
there is an integer $m$ such that $|M:K_i|=m^{\ell_i}$ for
$i=1,\ 2$. For any odd prime
$p$ we have $|M:K_1|_p=|M:K_2|_p$, which implies that $\ell_1=\ell_2$
and so $|K_1|=|K_2|$: a contradiction.

Suppose now that $T\cong\pomegap 82$. By Theorem~\ref{fullth}
$\sigma_1(K_i)'\times\cdots\times \sigma_k(K_i)'=K_i'$. We read off from
Table~\ref{full_simple} that in every case $|K_i:K_i'|$ is a 2-power, and
$|T_i:\sigma_i(K_j)|_5=1$. Therefore $|M:K_i|_5=k$ for $i=1,\ 2$, and so
$\ell_1=\ell_2$. This forces $|A|=|B|$ and inspection of
Table~\ref{full_simple} yields that $A\cong B\cong \sp
62\cong\Omega_7(2)$. And so line~3 of Table~\ref{groupsinttable} holds
with $q=2$.

Suppose that one of rows 1--3 of Table~\ref{groupsinttable} is
valid. The groups $A$ and $B$ in these rows are perfect, and so we
only have to show that $\cent{\sym\Omega}M=1$. 
By~\cite[Theorem~4.2A]{dm}, 
\begin{equation}\label{centM}
\cent {\sym \Omega}M\cong\norm
M{M_\omega}/M_\omega=\norm M{K_1\cap K_2}/(K_1\cap K_2).
\end{equation} 
It follows, however, from
Proposition~\ref{normcor1} that in this case $K_1\cap K_2$ is
self-normalising in $M$, and so $M$ is the unique minimal normal
subgroup of $G$. Thus $G$ is quasiprimitive. 
Suppose now that row~4 of Table~\ref{groupsinttable} is
valid. Then~\eqref{commnorm} follows from Theorem~\ref{fullth} and
Proposition~\ref{normcor1}. 
By Proposition~\ref{normcor1}
$$
\norm
M{K_1\cap K_2}=\norm M{K_1}\cap\norm M{K_2}.
$$ 
As $(\norm M{K_1}\cap\norm M{K_2})/(K_1'\cap K_2')$ is an elementary
abelian group of order $2^k$, by~\eqref{centM},  so is
$\cent{\sym\Omega}M$, and so all minimal normal subgroups of $G$
different from $M$ are also elementary abelian groups of order at most
$2^k$ .
\end{proof}

Finally in this section we show how to construct examples.

\begin{example}\label{fullex}
Let $T$ be a finite simple group with a non-trivial factorisation
$T=AB$, where $T$, $A$, and $B$ are as in Table~\ref{groupsinttable}. 
Set $\bar K_1=A^k$ and $\bar K_2=B^k$. Identify $M$ with $\inn M$ in
$\aut M$, and let 
$$
\bar G=M\left(\norm {\aut M}{\bar K_1}\cap \norm {\aut M}{\bar K_2}\right).
$$
Since the cyclic subgroup of $\aut M$ generated by the 
automorphism 
$$
\tau:(x_1,\ldots,x_k)\mapsto(x_k,x_1,\ldots,x_{k-1})
$$ is transitive on the set of simple direct factors of $M$ and
normalises $\bar K_1$, $\bar K_2$, we have that $M$ is a minimal normal subgroup of $\bar G$. Moreover,
since $\cent{\aut M}M=1$, we have that $M$ is the unique minimal normal
subgroup of $\bar G$. 

If $G_0=\norm{\aut M}{\bar K_1}\cap\norm{\aut M}{\bar K_2}$ then
$MG_0=\bar G$. 
As $\bar K_1$
and $\bar K_2$ are self-normalising in $M$, $M\cap G_0=\bar
K_1\cap \bar K_2$. 
Therefore, by~\cite[Lemma~4.1]{3types} the $M$-action on the coset space
$[M:\bar K_1\cap \bar K_2]$ can be extended to $\bar G$ with point stabiliser
$G_0$.
Moreover, $\bar K_1$ and $\bar K_2$ form a Cartesian system
for $M$ acted upon trivially by $G_0$. Consequently this action of $\bar
G$ preserves an intransitive $\bar G$-invariant Cartesian decomposition
given by the Cartesian system
$\{\bar K_1,\bar K_2\}$.
\end{example}

\section{The proof of Theorem~\ref{main}}\label{proof}

In this section we prove Theorem~\ref{main} working 
with the notation introduced in Section~\ref{mainthsec}. 

\begin{lemma}\label{max1}
Let $T_1,\ldots,T_k$, $K_1,\ldots,K_s$, $\bar\Xi_1,\ldots,\bar\Xi_s$,
and $\Omega_1,\ldots,\Omega_s$ be as in Section~$\ref{mainthsec}$. If, for some $i\in\{1,\ldots,k\}$ and $j\in\{1,\ldots,s\}$, 
$\sigma_i(K_j)$  is a proper maximal subgroup of $T_i$, then 
$\bar\Xi_j\in\cd{1}{G^{\Omega_j}}$. 
\end{lemma}
\begin{proof}
The group $G_\omega$ is transitive by conjugation on the set
$\{T_1,\ldots,T_k\}$, and, by Proposition~\ref{1ststep}, each of the $K_j$ is normalised by
$G_\omega$. Thus it suffices to prove that if 
$\sigma_1(K_1)$  is a proper maximal subgroup of $T_1$, then 
$\bar\Xi_1\in\cd{1}{G^{\Omega_1}}$. Assume without loss of generality
that $\Xi_1=\{\Gamma_1,\ldots,\Gamma_m\}$. 
By Proposition~\ref{1ststep}, the
$G$-action on $\Omega_1$ is equivalent to the $G$-action on $[M:K_1]$,
and $\bar\Xi_1\in\cd{\rm tr}{G^{\Omega_1}}$. Thus if
$\bar\Xi_1\in\cd{\rm S}{G^{\Omega_1}}$ then $\sigma_1(K_1)=T_1$. If 
$\bar\Xi_1\in\cd{2\sim}{G^{\Omega_1}}\cup\cd{2\not\sim}{G^{\Omega_1}}\cup\cd{3}{G^{\Omega_1}}$
then there are distinct $j_1,\ j_2\in\{1,\ldots,m\}$ such that
$\sigma_1(L_{j_1}),\ \sigma_1(L_{j_2})<T_1$, 
$\sigma_1(L_{j_1})\sigma_1(L_{j_2})=T_1$, and $\sigma_1(K_1)\leq
\sigma_1(L_{j_1})\cap \sigma_1(L_{j_2})$. Hence $\sigma_1(K_1)$ is not
a maximal subgroup of $T_1$. 

Suppose finally that $\bar\Xi_1\in\cd{\rm
1S}{G^{\Omega_1}}$. Then,
by~\cite[Theorem~6.1]{transcs}, we may assume without 
loss of generality that there  is a full strip $X$ of length~$2$ involved in $L_1$ 
covering $T_1$ and $T_2$, and there are indices $j_1,\ j_2\in\{2,\ldots,m\}$ 
such that $\sigma_1(L_{j_1})<T_1$, $\sigma_2(L_{j_2})<T_2$. Let $\alpha:
T_1\rightarrow T_2$ be the isomorphism such that 
$X=\{(t,\alpha(t))\ |\ t\in T_1\}$. 
It follows from~\cite[Lemma~2.1]{charfact} that 
$\sigma_1(L_{j_1})\alpha^{-1}(\sigma_2(L_{j_2}))=T_1$. In particular
$\sigma_1(L_{j_1})$ and $\alpha^{-1}(\sigma_2(L_{j_2}))$ are distinct
subgroups of $T_1$. On the other hand,
$\sigma_1(K_1)\leq\sigma_1(L_1\cap L_{j_1}\cap L_{j_2})\leq
\sigma_1(L_{j_1})\cap \alpha^{-1}(\sigma_2(L_{j_2}))$. Hence
$\sigma_1(K_1)$ cannot be a maximal subgroup of $T_1$. Therefore the
only remaining possibility is that 
$\bar\Xi_1\in\cd{1}{G^{\Omega_1}}$.
\end{proof}

Recall
that $\{L_1,\ldots,L_{\ell}\}$ is the original Cartesian system
corresponding to the intransitive Cartesian decomposition $\E$. The
following lemma is an easy consequence of~\cite[Lemma~3.1]{recog}.

\begin{lemma}\label{mixLK}
Let $L_1,\ldots,L_\ell$ be as in Theorem~$\ref{main}$, 
and suppose that $I_1,\ldots,I_m$ are pairwise disjoint subsets of
$\{1,\ldots,\ell\}$, and, for $i=1,\ldots,m$, set $Q_i=\bigcap_{j\in
I_i} L_j$. Then 
$$
Q_i\left(\bigcap_{j\neq
i}Q_j\right)=M\quad \mbox{for all}\quad i\in\{1,\ldots,m\}.
$$
\end{lemma}

\begin{proof}[The proof of Theorem~$\ref{main}$]
By Proposition~\ref{1ststep}, the partitions
$\Omega_i$ are $G$-invariant and $\bar\E=\{\Omega_1,\ldots,\Omega_s\}$
is a Cartesian decomposition of $\Omega$ on which $G$ acts trivially.
By the same result, 
for $i=1,\ldots,s$, the subgroup $K_i$ is the stabiliser in $M$ of the 
block in $\Omega_i$ containing $\omega$,  $\bar\Xi_i\in\cd{\rm
  tr}{G^{\Omega_i}}$, and $M$ is faithful on $\Omega_i$. 
It follows from
Theorem~\ref{3th} that the number $s$ of $G$-orbits on $\E$ is at
most~3. 

We prove the rest of Theorem~\ref{main} part by part.

(i) Suppose first, without loss of generality, that 
$\bar\Xi_1\in\cd{\rm S}{G^{\Omega_1}}$. Then by Theorem~\ref{5class}(a), 
$K_1$ is a subdirect subgroup of $M$ and 
it follows from Theorem~\ref{intrstripd2} that
$s=2$, and that $(M,K_1,K_2)$ is a full
strip factorisation. In particular, $\sigma_i(K_2)$ is a
maximal subgroup of $T_i$ for all $i$, and  
hence Lemma~\ref{max1} implies
that $\bar\Xi_2\in\cd 1{G^{\Omega_2}}$, as required.

(ii) Next assume without loss of generality 
that $\bar\Xi_1\in \cd{2\not\sim}{G^{\Omega_1}}$, and that
$\Xi_1=\{\Gamma_1,\ldots,\Gamma_m\}$. Note that there are $j_1,\
j_2\in\{1,\ldots,m\}$ such that $\sigma_1(K_{j_1}),\
\sigma_1(K_{j_2})<T_1$. 
If $\sigma_1(K_2)=T_1$, then, by Theorem~\ref{intrstripd2}, 
$\bar\Xi_2\in\cd{\rm S}{G^{\Omega_2}}$,
and so part~(i) implies that $\bar\Xi_1\in\cd 1{G^{\Omega_1}}$, which
is a contradiction. Hence $\sigma_1(K_2)<T_1$. 
If $s\geq 3$ then 
the same argument shows that $\sigma_1(K_3)<T_1$ and,
by Lemma~\ref{mixLK}, 
$\sigma_1(L_{j_1}),\ 
\sigma_1(L_{j_2}),\ \sigma_1(K_2),\ \sigma_1(K_3)$ form a strong multiple
factorisation of the finite simple group $T_1$. As, by
Theorem~\ref{th2}, 
such a factorisation has
at most 3~subgroups, 
this yields that $s=2$. Similarly, 
if there are two indices $j_3,\ j_4\in\{m+1,\ldots,\ell\}$ such
that $\sigma_1(L_{j_3}),\ \sigma_1(L_{j_4})<T_1$ then the subgroups
$\sigma_1(L_{j_1}),\ 
\sigma_1(L_{j_2}),\ \sigma_1(L_{j_3}),\ \sigma_1(L_{j_4})$ form a strong
multiple factorisation of $T_1$. 
This again is a contradiction and so $\bar\Xi_2\in\cd{1}{G^{\Omega_2}}\cup\cd{\rm 1S}{G^{\Omega_2}}$. 
Thus there is a unique index $j_3\in\{m+1,\ldots,\ell\}$ such that $\sigma_1(L_{j_3})<T_1$. 
The subgroups 
$\sigma_1(L_{j_1}),\ \sigma_1(L_{j_2}),\ \sigma_1(L_{j_3})$ form a strong
multiple factorisation of $T_1$ and so $T_1$ and these subgroups are as in 
Table~\ref{tablemult}. If 
$\bar\Xi_2\in\cd{\rm 1S}{G^{\Omega_2}}$ then, by Theorem~\ref{5class}(b), 
$T_1$ must also be as in 
Table~\ref{isomfact} and so $T_1\cong\pomegap 83$.
Further, $\sigma_1(L_{j_3})$, and hence $\sigma_1(K_2)$, 
must be a  maximal subgroup of $T_1$. This,  however, cannot be
the case if $\bar\Xi_2\in\cd{\rm 1S}{G^{\Omega_2}}$, by Lemma~\ref{max1}. 
Thus the assertions in part~(ii) all hold. 

(iii) Suppose without loss of generality that $\bar\Xi_1\in\cd{\rm 1S}{G^{\Omega_1}}\cup\cd{2\sim}{G^{\Omega_1}}\cup\cd 3{G^{\Omega_1}}$ 
and that $\Xi_1=\{\Gamma_1,\ldots,\Gamma_m\}$. It follows from
part~(i) that $\bar\Xi_i\not\in\cd{\rm 1S}{G^{\Omega_i}}$ for all
$i\in\{2,\ldots,s\}$. Thus for $i\in\{1,\ldots,k\}$ and
$j\in\{1,\ldots,s\}$ the projection $\sigma_i(K_j)$ is proper in
$T_i$. If $\bar\Xi_1\in\cd 3{G^{\Omega_1}}$ then there are pairwise
distinct indices 
$j_1,\ j_2,\ j_3\in\{1,\ldots,m\}$ such that $\sigma_1(L_{j_1}),\
\sigma_1(L_{j_2}),\ \sigma_1(L_{j_3})<T_1$. By Lemma~\ref{mixLK}, the subgroups 
$\sigma_1(L_{j_1}),\ \sigma_1(L_{j_2}),\ \sigma_1(L_{j_3}),\
\sigma_1(K_2)$ form a strong 
multiple factorisation of $T_1$, which is a contradiction, by
Theorem~\ref{th2}. Thus $\bar\Xi_1\not\in\cd 3{G^{\Omega_1}}$. 

Suppose next that $\bar\Xi_1\in \cd{2\sim}{G^{\Omega_1}}$. 
Then there are
distinct 
indices $j_1,\ j_2\in\{1,\ldots,m\}$ such that $\sigma_1(L_{j_1})$ and
$\sigma_1(L_{j_2})$ are proper isomorphic subgroups of $T_1$.
On the other hand, as $\sigma_1(K_2)<T_1$, 
the subgroups 
$\sigma_1(L_{j_1}),\ \sigma_1(L_{j_2}),\ \sigma_1(K_2)$ 
form a strong multiple factorisation of $T_1$. By Table~\ref{tablemult} 
such a factorisation cannot contain two isomorphic subgroups, 
and so this is a contradiction. Thus $\bar\Xi_1$
cannot be an element of $\cd {2\sim}{G^{\Omega_1}}$.

Suppose finally that $\bar\Xi_1\in \cd{\rm 1S}{G^{\Omega_1}}$. Then,
by~\cite[Theorem~6.1]{transcs}, we may assume without 
loss of generality that there  is a full strip $X$ of length~$2$ involved in $L_1$ 
covering $T_1$ and $T_2$, and there are indices $j_1,\ j_2\in\{2,\ldots,m\}$ 
such that $\sigma_1(L_{j_1})<T_1$, $\sigma_2(L_{j_2})<T_2$. Let $\alpha:
T_1\rightarrow T_2$ be the isomorphism such that 
$X=\{(t,\alpha(t))\ |\ t\in T_1\}$. It follows from~\cite[Lemma~2.1]{charfact} that 
$\sigma_1(L_{j_1})\alpha^{-1}(\sigma_2(L_{j_2}))=T_1$. 
Theorem~\ref{intrstripd2} and part~(i) implies
that $\sigma_1(K_2)<T_1$. As $(L_1\cap L_{j_1}\cap L_{j_2})K_2=M$ and
$$
\sigma_1(L_1\cap L_{j_1}\cap L_{j_2})\leq \sigma_1(L_{j_1})\cap
\alpha^{-1}(\sigma_2(L_{j_2})), 
$$
we obtain that
$(\sigma_1(L_{j_1})\cap\alpha^{-1}(\sigma_2(L_{j_2})))\sigma_1(K_2)=T_1$. 
Then~\cite[Lemma~4.3(iii)]{bad:fact} implies  that $(T_1,\{\sigma_1(L_{j_1}),\alpha^{-1}(\sigma_2(L_{j_2})),\sigma_1(K_2)\})$ is a strong
multiple factorisation. By Table~\ref{tablemult} 
distinct subgroups in such
a factorisation cannot be isomorphic. This is a contradiction, and so $\bar\Xi_1\not\in\cd {\rm 1S}{G^{\Omega_1}}$.

(iv) Suppose that $\E$ is homogeneous. Then it follows from Theorem~\ref{homth}
that $G$ has exactly 2~orbits on $\E$ and so $s=2$. The same result
implies that $K_1,\ K_2$ form a full factorisation of
$M$, and that $\sigma_i(K_j)$ is a maximal
subgroup of $T_i$, for each $i$ and $j$. 
Thus Lemma~\ref{max1} gives $\bar\Xi_i\in\cd 1{G^{\Omega_i}}$ for
$i=1,\ 2$.

(v) Finally suppose that $s=3$. 
By part~(i) $\bar\Xi_i\not\in\cd{\rm S}{G^{\Omega_i}}$ and, by part~(iii),
$\bar\Xi_i\not\in\cd{\rm 1S}{G^{\Omega_i}}$ for $i=1,\ 2,\ 3$.
If $\bar\Xi_i\not\in\cd{1}{G^{\Omega_i}}$ for some $i$ then there must
be 4~pairwise distinct indices $j_1,\ j_2,\ j_3,\
j_4\in\{1,\ldots,\ell\}$ such that $\sigma_1(L_{j_1})$,
$\sigma_1(L_{j_2})$, $\sigma_1(L_{j_3})$, $\sigma_1(L_{j_4})$ are
proper subgroups of $T_1$. By~\eqref{simpfact}, these subgroups form a
strong multiple factorisation of $T_1$, which is a contradiction, by
Theorem~\ref{th2}. Thus $\bar\Xi_i\in\cd 1{G^{\Omega_i}}$ for $i=1,\
2,\ 3$. It also follows from Theorem~\ref{3th} that $(M,\{K_1,K_2,K_3\})$
is a strong multiple factorisation.
\end{proof}

\end{document}